\documentclass[a4paper,12pts]{amsart}
\usepackage{graphicx} 

\usepackage{amsthm}
\usepackage{amssymb}
\usepackage{amsfonts}
\usepackage{mathtools}
\usepackage{mathrsfs}
\usepackage{graphicx}
\usepackage{color}
\usepackage{enumitem}
\usepackage{cancel}
\usepackage{setspace}
\usepackage[english]{babel}
\usepackage{csquotes}
\usepackage{etoolbox}

\usepackage{thmtools}
\usepackage{yfonts}
\usepackage{mfirstuc}

\usepackage{fancyvrb}

\usepackage{lipsum}

\newcommand\blfootnote[1]{%
  \begingroup
  \renewcommand\thefootnote{}\footnote{#1}%
  \addtocounter{footnote}{-1}%
  \endgroup
}

\usepackage{textgreek} 

\usepackage{amsmath}
\usepackage{amsthm, amssymb, amsfonts}
\usepackage{graphicx, color}
\usepackage{caption, subcaption} 
\usepackage{multicol} 
\usepackage{stmaryrd}
\usepackage{tikz-cd} 
\usepackage{lipsum}  
\usepackage[capitalise,noabbrev, nameinlink]{cleveref} 
\Crefname{part}{Part}{Parts}
\Crefname{step}{Step}{Steps}
\Crefname{prop}{Proposition}{Propositions}
\Crefname{prob}{Problem}{Problems}

\newcommand{\Fin}{\text{Fin}}

\newcommand{\mA}{\mathcal{A}}

\newcommand{\mL}{\mathcal{L}}
\newcommand{\mR}{\mathcal{R}}
\newcommand{\mE}{\mathcal{E}}
\newcommand{\mD}{\mathcal{D}}
\newcommand{\gap}[3]{({#1}_{#3},{#2}_{#3})_{\alpha\in\omega_1}}

\usepackage{enumitem}
\renewcommand{\leq}{\leqslant}
\renewcommand{\geq}{\geqslant}

\renewcommand{\epsilon}{\varepsilon}

\newcommand{\stickT}{%
\setbox255=\hbox{\raise1ex\hbox{$\hspace{0.2pt}\,\bullet\,$}}
\mathord{\rlap{\hbox to\wd255{\hss\hbox{$|$}\hss}}
\box255}
}
\newcommand{\stickS}{%
\setbox255=\hbox{\raise0.6ex\hbox{$\scriptstyle\bullet$}}
\mathord{\rlap{\hbox to\wd255{\hss\hbox{$\scriptstyle|$}\hss}}
\box255}
}

\numberwithin{equation}{section}
\theoremstyle{plain}
\newtheorem{theorem}[equation]{Theorem}
\newtheorem{lemma}[equation]{Lemma}
\newtheorem{proposition}[equation]{Proposition}
\newtheorem{corollary}[equation]{Corollary}

\theoremstyle{definition}
\newtheorem{definition}[equation]{Definition}

\newtheorem{problem}[equation]{Problem}

\theoremstyle{remark}

\newtheorem{rem}[equation]{Remark}

\usepackage{mathtools}
\usepackage{amsfonts}
\usepackage{amssymb}
\usepackage{amsthm}
\usepackage{amsmath}
\usepackage{dirtytalk}
\usepackage{mathrsfs}  
\usepackage{lineno}

 \usepackage{multicol}  
 \usepackage{bbm}
 \usepackage{tabularx}
\newenvironment{claimproof}[1][\proofname]
{\begin{proof}[#1]}
{\end{proof}}

\usepackage{xargs}                      
\usepackage[colorinlistoftodos,prependcaption,textsize=tiny]{todonotes}
\newcommandx{\unsure}[2][1=]{\todo[linecolor=red,backgroundcolor=red!25,bordercolor=red,#1]{#2}}
\newcommandx{\change}[2][1=]{\todo[linecolor=blue,backgroundcolor=blue!25,bordercolor=blue,#1]{#2}}
\newcommandx{\info}[2][1=]{\todo[linecolor=OliveGreen,backgroundcolor=OliveGreen!25,bordercolor=OliveGreen,#1]{#2}}
\newcommandx{\improvement}[2][1=]{\todo[linecolor=Plum,backgroundcolor=Plum!25,bordercolor=Plum,#1]{#2}}
\newcommandx{\thiswillnotshow}[2][1=]{\todo[disable,#1]{#2}}
\title{Forcing axioms and construction schemes }
\author{Jorge Cruz, Osvaldo Guzm\'an and Stevo Todor\v{c}evi\'c}
\begin{document}
\maketitle\begin{abstract}We continue the development of the theory of construction schemes over $\omega_1$ as introduced by the third author by studying their relation with forcing axioms. Formally, we introduce the cardinals $\mathfrak{m}^n_{\mathcal{F}}$ and use the consistency of $\mathfrak{m}^2_\mathcal{F}>\omega_1$ to prove a fundamental result relating gaps and almost disjoint families over $\omega$. The cardinals $\mathfrak{m}_\mathcal{F}$ are also used to prove some limiting results for contstruction schemes, some of which answer questions from \cite{schemescruz}. Finally, we show that PID implies the non-existence of $2$-capturing schemes.
\end{abstract}

 \blfootnote{The first and second authors were partially supported by a PAPIIT grant IN101323 and CONACyT grant A1-S-16164. The second author was also supported by PAPIIT grant IA102222. The third author is partially supported by grants from NSERC (455916), CNRS (IMJ-PRG-UMR7586) and SFRS (7750027-SMART).}Experience has shown that \emph{forcing axioms} (such as MA, PFA, MM, \dots ) and guessing principles (for example $\diamond$, $\diamond^+$, $\diamond^*$,\dots) frequently have diametrically opposed consequences, thus giving rise to completely distinct universes of set theory. There is a plethora of mathematical statements $\phi$ whose independence from ZFC was proved by realizing that $\phi$ can be deduced from a certain forcing axiom and $\neg\phi$ can be deduced from a particular guessing principle. A famous example of this is the independence of \emph{Whitehead's problem} due to Shelah (see \cite{infiniteabeliangroups}). Evidently, this strategy has its limitations. However,  whenever we suspect a certain statement to be independent, it usually a good idea to ask about its relation with these types of axioms. Nevertheless, there are independent statements for which the only consistency proofs we know of, rely on combining components of both guessing principles and forcing axioms. An example of this phenomenon is the famous \emph{Katetov's metrization problem}, which asks whether every compact space whose square is $T_5$ (hereditarily normal) is in fact metrizable. The (consistent) positive answer to this problem was provided by Larson and the third author in \cite{katetovproblem}. There they introduce the forcing axiom MA$_{\omega_1}(T)$ (for a Suslin tree $T$), which can be thought as the largest protion of Martin's axiom which is compatible with the existence of $T$  (see also \cite{pmaxvariation} and \cite{chainconditionsmaximalmodels}). The  consistency of Katetov's problem is then obtained by forcing with $T$ over a model of MA$_{\omega_1}(T)$.
Since \cite{katetovproblem}, forcing axioms compatible with the existence of a Suslin tree have been widely used and studied (see \cite{pfasdow}, \cite{pfasdowtall}, \cite{pfaslocallycompact}, \cite{pengpfas}, \cite{Cesaruniformization}, \cite{forcingcoherentsuslinlocallycountable}, \cite{rudindowkersuslin}, \cite{forcingaxiomnonspecial}). In this way, we consider that the study of maximal versions of forcing axioms compatible with guessing principles is very promising. In this paper, we present the maximal version of Martin's axiom compatible with \emph{construction schemes.}

In \cite{schemenonseparablestructures}, the third author introduced a generalization of Jensen's $(\omega,1)$-gap morasses which he called \emph{construction schemes}. These objects are special families of finite subsets of $\omega_1$ which serve as a tool for building uncountable objects by means of finite approximations. In that same paper, he studied construction schemes with certain guessing properties which can be interpreted as finite dimensional versions of the properties possessed by a $\Diamond$-sequence. These type of schemes are called fully capturing, capturing and $n$-capturing (here, $n$ varies over $\omega$). It is easy to prove  (see \cite{schemescruz}) that the existence of $n$-capturing construction schemes is incompatible with MA.
In this paper, given an $n$-capturing construction scheme $\mathcal{F}$, we will define its $n$-parameterized Martin's number $\mathfrak{m}^n_\mathcal{F}$ as the Martin's number corresponding to the family of $ccc$-forcing notions which force our construction scheme to be $n$-capturing. We will show that the principle MA$^n(\mathcal{F})$ (stating that $\mathfrak{m}^n_\mathcal{F}=\mathfrak{c}>\omega_1$) is consistent with ZFC. In order to do this, we will use a finite support iteration of $ccc$-forcings. This approach will require to prove iteration theorems regarding the preservation of $n$-capturing construction schemes. As a main application of this principle, we will show the independence of a statement related to gaps over $\omega$ which not only PFA and CH decide in decide in the negative, but also axioms such as  PFA$(T)$ and PFA$(T)[T]$(we will deal with these two axioms in an upcoming paper). Before we continue the discussion, it will be convenient to remind  some important notions.

Let $\mL$ and $\mR$ be two families of subsets of a set $X$. The pair $(\mL,\mR)$ is said to be a \emph{pregap} if $L\cap R=^*\emptyset$ for all $L\in \mL$ and $R\in \mR$. A subset $C$ of $X$ is said to \textit{separate} $(\mL,\mR)$ if\footnote{Here $L\subseteq^*C$ and  $C\cap R=^*\emptyset $ mean that $L\backslash C$ and $C\cap R$ are finite.} $L\subseteq^*C$ and $C\cap R=^*\emptyset$ for each $L\in\mL$ and $R\in \mR$. We say that $(\mL,\mR)$ is a \emph{gap} in case it is a pregap and there is no $C\subseteq X$ separating it. Of particular interest for us are the pregaps (gaps) $(\mL,\mR)$ for which both $\mL$ and $\mR$ are towers, that is, they are well ordered with respect to $\subseteq^*$. If this situation occurs and $\kappa$ and $\lambda$ are the order types of $\mathcal{L}$ and $\mathcal{R}$ respectively, we refer to $(\mL,\mR)$ as a $(\kappa,\lambda)$-pregap (gap). Whenever $(\mathcal{L},\mathcal{R})$ is a pregap and $\mathcal{L}$ and $\mathcal{R}$ are enumerated as $\langle L_\alpha\rangle_{\alpha\in\kappa}$ and $\langle R_\alpha\rangle_{\alpha\in\kappa}$ respectively, we  rewrite $(\mathcal{L},\mathcal{R})$ as $(L_\alpha,R_\alpha)_{\alpha\in \kappa}$. Moreover, if $(\mathcal{L},\mathcal{R})$ is assumed to be a $(\kappa,\kappa)$-pregap, we always assume such enumerations to be increasing with respect to $\subseteq^*$.\\
 Recall  that an \emph{almost disjoint family} (shortly, AD family) over a set $X$,\footnote{If $X $ is not mentioned, we will assume that it is $\omega$.} is a family $\mathcal{A}$ of infinite subsets of $X$ so that $A\cap B=^*\emptyset$ for any two distinct $A,B\in \mathcal{A}$. We say  $\mathcal{A}$ is \emph{inseparable} if for any  two uncountable pairwise disjoint $\mD,\mE\subseteq  \mathcal{A}$, the pair $(\mD,\mE)$ forms a gap\footnote{ In \cite{luzinandantiluzin}, Soukup and Roitman call these families \emph{almost Luzin.}}. AD families are one of the central objects of study in modern combinatorial set theory.  Constructing  almost disjoint families with special properties is usually difficult and in most cases had lead to the development of powerful tools (see  \cite{InvariancePropertiesofAlmostDisjointFamilies}, \cite{Madfamilieswithstrongcombinatorialproperties}, \cite{madnessandnormality},
\cite{frechetlike},  \cite{nonpartitionable}, \cite{cechfunction},\cite{ThereisavanDouwenMADfamily} and \cite{SANEPlayer}). Almost disjoint families have also played a central roll in the solution of many problems of Topology and Analysis. An example of this is the solution of the selection problem posed by van Mill and Wattel in \cite{van1981selections} and solved by  Hru\v{s}\'ak and Mart\'inez-Ru\'iz in \cite{hruvsak2009selections}. The reader interested in learning more about AD families is referred to \cite{TopologyofMrowkaIsbellSpaces}, \cite{Combinatoricsoffiltersandideals} and \cite{AlmostDisjointFamiliesandTopology}.\\
One of the main motivations of this work is to understand up to which degree do almost disjoint families determine the structure of gaps over $\omega$. In \cite{guzman2019mathbb}, the second author, Hru\v{s}ak and Koszmider constructed a special kind of inseparable AD family $\mathcal{A}$ fo size $\omega_1$. 
Among other things, they used $\mA$ to recursively construct an $(\omega_1,\omega_1)$-pregap $(\mL,\mR)=\gap{L}{R}{\alpha}$ with the key feature that  $L_{\alpha+1}\backslash L_\alpha, R_{\alpha+1}\backslash R_\alpha\in \mathcal{A}$ for any $\alpha\in \omega_1$. This automatically implies that $(\mL,\mR)$ is a gap. Indeed, any set separating $(\mL,\mR)$ would also separate $(L_{\alpha+1}\backslash L_\alpha, R_{\alpha+1}\backslash R_\alpha)_{\alpha\in \omega_1}$ which is impossible because $\mathcal{A}$ is inseparable. This provides us with a precise way of saying that the \say{gapness} of $(\mathcal{L},\mathcal{R})$ is determined by an almost disjoint family.

\begin{definition}[Levelwise-inseparable gaps]We call an $(\omega_1,\omega_1)$-gap $(L_\alpha,R_\alpha)_{\alpha\in\omega_1}$ \emph{levelwise-inseparable} if $(L_{\alpha+1}\backslash L_\alpha, R_{\alpha+1}\backslash R_\alpha)_{\alpha\in\omega_1}$ forms a gap.\footnote{Note that $\{L_{\alpha+1}\backslash L_\alpha\,:\,\alpha\in \omega_1\}\cup\{L_{\alpha+1}\backslash R_\alpha\,:\,\alpha\in \omega_1\}$ forms an AD family.} If a gap is not levelwise-inseparable, we call it \emph{levelwise-separable}.
\end{definition}
The notion of levelwise-separability is too strong. As we will see, there are levelwise-separable gaps (see Theorem \ref{firsthausdorffdonutseparable}). However, when this notion is weaken the situation changes dramatically.  Recall that a \emph{subgap} $(\mathcal{L}',\mathcal{R}')$ of a gap $(\mathcal{L},\mathcal{R})$ is a gap for which $\mathcal{L}'\subseteq \mathcal{L}$ and $\mathcal{R}'\subseteq \mathcal{R}$.
\begin{definition}[Weakly levelwise-inseparable gaps]We call an $(\omega_1,\omega_1)$-gap \emph{weakly levelwise-inseparable} (shortly, WLI) if it has a levelwise-inseparable  subgap.
\end{definition}
In this paper we will show  that the statement \say{Every $(\omega_1,\omega_1)$-gap is WLI} is independent from $ZFC.$ This result is interesting as it tells us how important is the roll played by almost disjoint families in the structure of gaps over $\omega$. It turns out that both CH (hence, also stronger guessing principles such $\Diamond$) and PFA imply that every gap is weakly levelwise-inseparable. As we prove in an upcoming paper, the same holds when we consider either of the principles PFA(T) or PFA(T)[T] where $T$ is a Suslin tree.    In order to show that consistently there is gap which is not WLI, we will show that this statement is implied by the inequality $\mathfrak{m}^2_\mathcal{F}>\omega_1$. Lastly, we will show that MA does not decide this problem. In order to do this, we will prove that non-WLI gaps can not be destroyed by $ccc$ forcing notions.
\\
 Beyond the application mentioned above, we will use the cardinals $\mathfrak{m}^n_\mathcal{F}$ to solve some important open problems from  \cite{schemescruz}. Namely, we will show that:\begin{itemize}
    \item Under $\mathfrak{m}^2_\mathcal{F}>\omega_1$ every Aronszajn tree is special and every $(\omega_1,\omega_1)$-gap is indestructible. In particular, the existence of a $2$-capturing scheme is not strong enough to imply the existence of neither destructible gaps nor Suslin trees. It was already known that destrucible gaps exist if there is a $3$-capturing scheme (see \cite{treesgapsscheme}) and Suslin trees exist if there is fully capturing construction scheme (see \cite{schemescruz}).
    \item It is independent from ZFC+\say{There is an $n$-capturing construction scheme} that every $n$-capturing scheme is $n$-capturing with partitions.
\end{itemize}
Related to the proof of the second item above is the fact that under $\mathfrak{m}^n_\mathcal{F}>\omega_1$ there is an ultrafilter $\mathcal{U}(\mathcal{F})$ over $\omega_1$ which has a $\Sigma^1_1$ definition with the scheme $\mathcal{F}$ as a parameter. Interestingly, any non-trivial projection of this ultrafilter to $\omega$ is a Ramsey ultrafilter. A similar phenomenom was studied by the third author in  \cite{mapsontrees} and \cite{walksultrafilters}. There he defined the ultrafilter $\mathcal{U}(T)$ associated to an Aronszajn tree $T$ and showed that $\mathcal{U}(T)$ behaves in a similar way to $\mathcal{U}(\mathcal{F})$ under  $\mathfrak{m}>\omega_1.$
Finally, we will show that the $P$-ideal dichotomy implies the non-existence of $2$-capturing schemes. In an upcoming paper, we will prove the consistency the axiom PFA$(\mathcal{F})$ (whose definition is as expected), and use it to show, among other things, that the mapping reflection principle MRP is consistent with the existence of a $2$-capturing scheme and that $\mathcal{U}(\mathcal{F})$ has maximal Tukey type among all ultrafilters over $\omega_1$. The consistency of this new axiom will be achieved by using Neeman's iteration method. \\

The paper is structured as follows: In Section \ref{everygapiswli} we prove that under CH and PFA every gap is weakly levelwise-inseparable. No knowledge of construction schemes will be necessary here. In Section \ref{sectionschemes} we provide a quick review of capturing construction schemes and present the parametrized Martin's numbers $\mathfrak{m}_\mathcal{F}^n$.  In Section \ref{noteverygapiswli} we show that under $\mathfrak{m}_\mathcal{F}>\omega_1$ there is a gap which is not WLI. In Section \ref{fragmentmartin}, we study in a greater detail the the cardinals $\mathfrak{m}_\mathcal{F}^n$. In Section \ref{sectiontreesandgaps} we prove that $\mathfrak{m}^2_\mathcal{F}>\omega_1$ implies that every Aronszajn tree is special and that every gap is indestructible. In Section \ref{capturingwithpartitionssection}, we show that a construction scheme which is $n$-capturing  may not be $\mathcal{P}$-$n$-capturing for any non-trivial partition of $\omega$. This answers a question from \cite{schemescruz}. In Section \ref{sectionramseyultafilters}, we show that Ramsey ultafilters exist under $\mathfrak{m}_\mathcal{F}^n>\omega_1$. In Section \ref{sectionpidschemes}, we show that the P-ideal dichotomy implies the non-existence of $2$-capturing schemes. Finally, in Section \ref{sectionproblems} we pose some open problems.
\section{Notation}
The notation and terminology used is here mostly standard and it follows \cite{schemescruz}.  Given a set $X$ and a (possibly finite) cardinal $\kappa$, $[X]^\kappa$ denotes the family of all subsets of $X$ of cardinality $\kappa$. The sets $[X]^{<\kappa}$ and $[X]^{\leq \kappa}$ have the expected meanings. The family of all non-empty finite sets of $X$ is denoted by $\Fin(X)$. That is, $\Fin(X)=[X]^{<\omega}\backslash\{\emptyset\}$. $\mathscr{P}(X)$ denotes the power set of $X$. By $\text{Lim}$  we mean the  set of limit ordinals strictly smaller than $\omega_1$. Given sets $X$ and $Y$ two sets of ordinals, we write $X<Y$ whenever $\max(X)<\min(Y)$ or $X=\emptyset$.

For a set of ordinals $X$, we denote by $ot(X)$ its order type. We identify $X$ with the unique strictly increasing function $h:ot(X)\longrightarrow X$. In this way, $X(\alpha)=h(\alpha)$ denotes the $\alpha$ element of $X$ with respect to its increasing enumeration. Analogously, $X[A]=\{X(\alpha)\,:\,\alpha\in A\}$ for $A\subseteq ot(X)$.  A family $\mathcal{D}$ is called a \textit{$\Delta$-system} with root $R$ if $|\mathcal{D}|\geq 2$ and $X\cap Y= R$ whenever $X,Y\in \mathcal{D}$ are different.  If moreover, $R<X\backslash R$ for any and $X\backslash R<Y\backslash R$ or viseversa for any two $X,Y\in \mathcal{D}$, we call $\mathcal{D}$ a root-tail-tail $\Delta$-system.

Given a forcing notion $\mathbb{P}$, we denote by $\mathfrak{m}(\mathbb{P})$ its Martin's number, i.e., the minimal cardinal $\kappa$ for which there are $\kappa$-many dense sets over $\mathbb{P}$ which are not intersected by any filter. We say that $\mathbb{P}$ is $n$-Knaster (or has poperty $K_n$) if for each uncountable $\mathcal{A}\subseteq\mathbb{P}$ we can find $\mathcal{B}\in [\mathcal{A}]^{\omega_1}$  which is $n$-linked. That is, every $p\in [\mathcal{B}]^n$ is bounded below inside $\mathbb{P}$. $\mathbb{P}$ is said to have precaliber $\omega_1$ if every subset of $\mathcal{\mathbb{P}}$ of size $\omega_1$ has an uncountable centered subset. The Martin's number for $n$-Knaster forcings is defined as $\mathfrak{m}_{K_n}=\min(\mathfrak{m}(\mathbb{P})\,:\,\mathbb{P}\text{ is $n$-Knaster})$. Finally, given a cardinal $\kappa$, $\mathbb{C}_\kappa$ stands for the Cohen forcing for adding $\kappa$-many Cohen reals.

\section{Every gap is weakly levelwise-inseparable}\label{everygapiswli}
In this section we aim to prove that consistently every gap is weakly levelwise-inseparable under CH, PFA and after forcing with $\mathbb{C}_k$ over a model of CH. The two latter proofs are based on the fact that, under CH, $(\omega_1,\omega_1)$-gaps satisfy a stronger property than just being weakly levelwise-inseparable. Let $(\mathcal{L},\mathcal{R})=(L_\alpha,R_\alpha)_{\alpha\in\omega_1}$ be an $(\omega_1,\omega_1)$-gap. A set $X\in [\omega_1]^{\omega_1}$ is said to be \emph{adequate} for $(\mathcal{L},\mathcal{R})$ if for any $Y\in[\omega_1]^{\omega_1}$, the pregap $$(L_{X(\alpha+1)}\backslash L_{X(\alpha)},R_{X(\alpha+1)}\backslash R_{X(\alpha)})_{\alpha\in Y}$$ is a gap. Note that if an $(\omega_1,\omega_1)$-pregap admits an adequate set, then such pregap is weakly levelwise-inseparable.

\begin{theorem}[Under CH]\label{weaklydonutch} Let $(\mathcal{L},\mathcal{R})=(L_\alpha,R_\alpha)_{\alpha\in\omega_1}$ be an $(\omega_1,\omega_1)$-gap. Then there is a club $X\in [\omega_1]^{\omega_1}$ adequate for $(\mathcal{L},\mathcal{R}).$ In particular, every gap is weakly levelwise-inseparable.
\begin{proof}Since we are assuming CH we can enumerate $[\omega]^{\omega}$ as $\langle C_\alpha\rangle_{\alpha\in \omega_1}$. We will build $X$ be recursion in such way that for any $\beta\in \omega_1$ and each $\alpha\leq \beta$ one of the following conditions occur:
\begin{center}\begin{minipage}{5cm} \begin{center} \textbf{(A)}\end{center} $$L_{X(\beta+1)}\backslash L_{X(\beta)}\not\subseteq^* C_\alpha.$$
\end{minipage}\hspace{2cm} \begin{minipage}{5cm}\begin{center} \textbf{(B)}\end{center} $$(R_{X(\beta+1)}\backslash R_{X(\beta)})\cap C_\alpha\not=^*\emptyset.$$
\end{minipage}
\end{center}
If $\gamma$ is limit and we have constructed $X(\alpha)$ for any $\alpha<\gamma$ we just define $X(\gamma)$ as $\sup(X(\alpha)\,:\,\alpha<\gamma)$. The interesting case is when we have constucted $X(\beta)$ for some $\beta\in \omega_1$ and we want to define $X(\beta+1).$ Here, we have the following claim:\\\\
\underline{Claim}: The pregap $(L_\delta\backslash L_{X(\beta)}, R_\delta\backslash R_{X(\beta)})_{\delta>X(\beta)}$ is a gap.
\begin{claimproof}[Proof of claim] Suppose towards a contradiction that this is not the case and let $C$ be a set separating $(L_\delta\backslash L_{X(\beta)}, R_\delta\backslash R_{X(\beta)})_{\delta>X(\beta)}$. Then $C\cup L_{X(\beta)}$ separates $(\mathcal{L},\mathcal{R})$ which is a contradiction. 
\end{claimproof}
Fix $\alpha\leq \beta$. According to the previous claim, the set $C_\alpha$ does not separate   $(L_\delta\backslash L_{X(\beta)}, R_\delta\backslash R_{X(\beta)})_{\delta>X(\beta)}$. Thus, we can find $\delta_\alpha >X(\beta)$ for which either $L_{\delta_\alpha}\backslash L_{X(\beta)}\not\subseteq^*C_\alpha$ or  $R_{\delta_\alpha}\cap C_\alpha\not=^*\emptyset.$ Let us define $X(\beta+1)$ as $\sup(\delta_\alpha\,:\alpha\leq \beta)$. Then $L_{\delta_\alpha}\backslash L_{X(\beta)}\subseteq L_{X(\beta+1)}\backslash L_{X(\beta)}$ and $R_{\delta_\alpha}\backslash R_{X(\beta)}\subseteq R_{X(\beta+1)}\backslash R_{X(\beta)}$ for each $\alpha\leq \beta$. In this way we guarantee that either condition (A) or condition (B) will occur for any such $\alpha.$ This finishes the recursion.

We will now prove that $X$ is as desired. Let $Y\in [\omega_1]^{\omega_1}$. We need to show that $(L_{X(\alpha+1)}\backslash L_{X(\alpha)},R_{X(\alpha+1)}\backslash R_{X(\alpha)})_{\alpha\in Y}$ is a gap.
For this, fix $C$ an infinite subset of $\omega$. By the assumptions, we know that there is $\alpha\in \omega_1$ so that $C=C_\alpha$. Let $\beta\in Y$ be such that $\beta>\alpha$. Therefore, either $L_{X(\beta+1)}\backslash L_{X(\beta)}\not\subseteq^*C_\alpha$ or $(R_{X(\beta+1)}\backslash R_{X(\beta)})\cap C_\alpha\not=^*\emptyset.$ In any case, $X(\beta)$ witness that $C_\alpha$ does not separate the gap that we are considering.
\end{proof}
\end{theorem}
Now we  analyze the relationship between levelwise-separable gaps and the Cohen forcing $\mathbb{C}$.
\begin{lemma}\label{adequatecohenpreserving} Let $(\mathcal{L},\mathcal{R})=(L_\alpha,R_\alpha)_{\alpha\in\omega_1}$ be an $(\omega_1,\omega_1)$-gap and let $X\in [\omega_1]^{\omega_1}$. If $X$ is adequate for $(\mathcal{L},\mathcal{R})$, then  $\mathbb{C}\Vdash \text{\say{ $X\textit{ is adequate for }(\mathcal{L},\mathcal{R})$}}.$
\begin{proof}Let $G$ be a $\mathbb{C}$-generic filter over $V$ and let $Y\in [\omega_1]^{\omega_1}\cap V[G]$. Suppose towards a contradiction that there is $C\in V[G]$ separating the pregap $$(L_{X(\alpha+1)}\backslash L_{X(\alpha)},R_{X(\alpha+1)}\backslash R_{X(\alpha)})_{\alpha\in Y}.$$ Since $Y$ is uncountable we can find $a,b\in [\omega]^{<\omega}$  and $Y'\in [Y]^{\omega_1}\cap V[G]$ so that for any $\alpha\in Y'$, $(L_{X(\alpha+1)}\backslash L_{X(\alpha)})\backslash C= a$ and $(R_{X(\alpha+1)}\backslash R_{X(\alpha)})\cap C=b$.
Due to well-known facts concerning the Cohen forcing $\mathbb{C}$, we know there is $Z\in [\omega_1]^{\omega_1}\cap V$ with $Z\subseteq Y'$. Let us define $$C'=\bigcup\limits_{\alpha\in Z}L_{X(\alpha+1)}\backslash L_{X(\alpha)}.$$
Then $C'\in V$ and trivially $L_{X(\alpha+1)}\backslash L_{X(\alpha)}\subseteq C'$ for any $\alpha\in Z$. Furthermore, $C'\subseteq C\cup a$. Therefore, $(R_{X(\alpha+1)}\backslash R_{X(\alpha)})\cap C'$ is finite for any $\alpha\in Z.$ We conclude that, in $V$, there is a set separating the gap $(L_{X(\alpha+1)}\backslash L_{X(\alpha)},R_{X(\alpha+1)}\backslash R_{X(\alpha)})_{\alpha\in Z}.$ This is a contradiction to the hypotheses. 
\end{proof}

\end{lemma}
\begin{theorem}[Under CH] Let $\kappa$ be an uncountable cardinal. Then $$\mathbb{C}_\kappa\Vdash\text{\say{ Every $(\omega_1,\omega_1)$-gap is weakly levelwise-inseparable }}.$$
In particular, the statement \say{Every $(\omega_1,\omega_1)$-gap is weakly levelwise-inseparable} is consistent with an arbitrarily large continuum.
\begin{proof}If $\kappa=\omega_1$ the argument is clear, since in this case $\mathbb{C}_\kappa$ forces CH. So suppose that $\kappa>\omega_1$ and let $G$ be a $\mathbb{C}_\kappa$-generic filter over $V$. Finally, let $(\mathcal{L},\mathcal{R})=(L_\alpha,R_\alpha)_{\alpha\in \omega_1}\in V[G]$ be an arbitrary $(\omega_1,\omega_1)$-gap. Since $|\mathcal{L}|=|\mathcal{R}|=\omega_1$, then there is a $\mathbb{C}_{\omega_1}$-generic filter over $V$, namely $H$, and  a $\mathbb{C}_{\kappa}$-generic filter over $V[H]$, namely $K$, so that $V[H][K]=V[G]$ and $(\mathcal{L},\mathcal{R})\in V[H].$ By the hypotheses, $V$ is a model of CH. Therefore $V[H]$ models CH too. According to Theorem \ref{weaklydonutch} there is $X\in V[H]$ which is adequate for $(\mathcal{L},\mathcal{R}).$ We finish by proving the following claim.\\\\
\underline{Claim}: $X$ testifies that $(\mathcal{L},\mathcal{R})$ is weakly levelwise-inseparable.
\begin{claimproof}[Proof of claim.] Suppose towards a contradiction that there is $C\in V[G]$ separating the pregap $$(L_{X(\alpha+1)}\backslash L_{X(\alpha)}, R_{X(\alpha+1)}\backslash R_{X(\alpha)})_{\alpha\in \omega_1}.$$
Again, since $|C|=\omega$ we can find $\mathbb{C}$-generic filter over $V[H]$, say $H'$, for which $C\in V[H][H']$. By virtue of Lemma \ref{adequatecohenpreserving}, $V[H][H']\models X\textit{ is adequate for }(\mathcal{L},\mathcal{R}).$ In particular $C$ can not separate the pregap that we are considering. This contradiction finishes the proof.
\end{claimproof}
\end{proof}
\end{theorem}

We will finish this section by showing  that PFA also implies that every gap is weakly levelwise-inseparable. For this sake, we need to generalize the concept of an adequate set. Let us say that a pregap $(\mathcal{D},\mathcal{E})=(D_\alpha,E_\alpha)_{\alpha\in\omega_1}$ (not necessarily of type $(\omega_1,\omega_1))$ is \emph{almost Luzin} if $(\mathcal{D}',\mathcal{E}')$ is a gap for any two uncountable $\mathcal{D}'\subseteq \mathcal{D}$ and $\mathcal{E}'\subseteq \mathcal{E}$.
\begin{definition}[Highly adequate sets]Let $(\mathcal{L},\mathcal{R})=(L_\alpha,R_\alpha)_{\alpha\in\omega_1}$ be an $(\omega_1,\omega_1)$-gap. We say that $X\in [\omega_1]^{\omega_1}$ is \textit{highly adequate} for $(\mathcal{L},\mathcal{R})$ if the pregap $$(L_{X(\alpha+1)}\backslash L_{X(\alpha)},R_{X(\alpha+1)}\backslash R_{X(\alpha)})_{\alpha\in X}$$
is almost Luzin.
\end{definition}
Trivially, if $X$ is highly adequate for $(\mathcal{L},\mathcal{R})$ then it is also adequate.  It is very easy to check that the dichotomy satisfied by the club $X$ that we constructed in Theorem \ref{weaklydonutch} already implies that  $X$ is a highly adequate set for the gap $(\mathcal{L},\mathcal{R})$. Hence, we have the following proposition.
\begin{proposition}[Under CH]\label{highlyadequateprop} Let $(\mathcal{L},\mathcal{R})=(L_\alpha,R_\alpha)_{\alpha\in\omega_1}$ be an $(\omega_1,\omega_1)$-gap. Then there is a club $X\in [\omega_1]^{\omega_1}$ highly adequate for $(\mathcal{L},\mathcal{R})$.
\end{proposition}

\begin{definition}[Normal pregaps]\label{normalgapdef} Suppose that $
(D_\alpha,E_\alpha)_{\alpha\in\omega_1}$ is a pregap. We say that it is \textit{normal} if $D_\alpha\cap E_\alpha=\emptyset$ for each $\alpha\in\omega_1.$
    
\end{definition}
Suppose that $
(D_\alpha,E_\alpha)_{\alpha\in\omega_1}$ is a pregap and consider, for each $\alpha$, $D'_\alpha=D_\alpha\backslash E_\alpha$ and $E'_\alpha=E_\alpha\backslash D_\alpha$. Then $(D'_\alpha,E'_\alpha)_{\alpha\in\omega_1}$ is a normal pregap. Furthermore, both gaps induced the same equivalence classes in $\mathscr{P}(\omega)/\Fin$. In particular, $
(D_\alpha,E_\alpha)_{\alpha\in\omega_1}$ is WLI if and only if the same is true for $
(D'_\alpha,E'_\alpha)_{\alpha\in\omega_1}$.

\begin{definition}[Biorthogonal gaps] Let $(\mathcal{D},\mathcal{E})=(D_\alpha,E_\alpha)_{\alpha\in \omega_1}$ be a normal pregap. We say that $(\mathcal{D},\mathcal{E})$ is \textit{biorthogonal} if $(D_\alpha\cap E_\beta)\cup (D_\beta\cap E_\alpha)\not=\emptyset$ for all $\alpha\not=\beta\in \omega_1$.
\end{definition}
The following is already known for $(\omega_1,\omega_1)$-gaps.
\begin{lemma}\label{biorthogonalgaplemma}If $(\mathcal{D},\mathcal{E})=(E_\alpha,D_\alpha)_{\alpha\in\omega_1}$ is a normal biorthogonal pregap, then it is a gap.
\begin{proof}Suppose towards a contradiction that there is $C\in [\omega]^{\omega}$ which separates $(\mathcal{D},\mathcal{E})$. According to the pigheonhole principle, there is $X\in [\omega_1]^{\omega_1}$ for which $D_\alpha\backslash C=D_\beta\backslash C$ and $E_\alpha\cap C=E_\beta\cap C$ for all $\alpha,\beta\in X.$ Fix two distinct $\alpha,\beta\in X$. Then \begin{align*}D_\alpha\cap E_\beta&= (D_\alpha\cap (E_\beta\cap C))\cup (E_\beta\cap (D_\alpha\backslash C))\\&=(D_\alpha\cap (E_\alpha\cap C))\cup (E_\beta\cap (D_\beta\backslash C))=\emptyset
\end{align*}
By symmetry, we also have that $D_\beta\cap E_\alpha=\emptyset$, but this is a contradiction to the normality of $(\mathcal{D},\mathcal{E})$. Thus, the proof is over.
\end{proof}
    
\end{lemma}

Let $(\mathcal{D},\mathcal{E})=(D_\alpha,E_\alpha)_{\alpha\in \omega_1}$ be a normal gap. We define the forcing $\mathbb{L}(\mathcal{D},\mathcal{E})$ as the set of all $p\in [\omega_1]^{<\omega}$ so that $(D_\alpha\cap E_\beta)\cup (D_\beta\cap E_\alpha)\not=\emptyset$
for all $\alpha\not=\beta\in p.$ The order is given by  $$p\leq q\text{ if and only if }q\subseteq p.$$

\begin{proposition}Let $(\mathcal{D},\mathcal{E})=(D_\alpha,E_\alpha)_{\alpha\in \omega_1}$ be a normal gap.  If $(\mathcal{D},\mathcal{E})$ is almost Luzin, then $\mathbb{L}(\mathcal{D},\mathcal{E})$ is $ccc$. 
\begin{proof}
Let $\mathcal{A}$ be an uncountable subset of $\mathbb{L}(\mathcal{D},\mathcal{E})$. We will show that $\mathcal{A}$ is not an antichain.  By refining $\mathcal{A}$, we may assume without loss of generality that there are $n,m\in \omega$ for which the following properties hold for each $p,q\in \mathcal{A}$:
\begin{enumerate}[label=$(\arabic*)$]
\item $|p|=n$.
\item For all $i,j<n$, $(D_{p(i)}\cap E_{p(j)})\cup (D_{p(j)}\cap E_{p(i)})\subseteq m$.
\item For all $i<n$, $D_{p(i)}\cap m=D_{q(i)}\cap m$ and $D_{p(i)}\cap m=D_{q(i)}\cap m$.\end{enumerate}
Note that for any two given conditions in $p,q\in \mathbb{L}(\mathcal{D},\mathcal{E})$, we have that $p$ and $q$ are compatible if and only if $p\backslash q$ and $q\backslash p$ are compatible conditions. Because of this and due to the $\Delta$-system Lemma, we may also assume that the elements of $\mathcal{A}$ are pairwise disjoint. Now, we proceed to find  distinct $p,q\in \mathcal{A}$ for which $p\cup q$ is a condition. This will be done after proving the following claims. \\

\noindent
\underline{Claim 1}: Let $p,q\in \mathcal{A}$ and $i\not=j<n$, then $(D_{p(i)}\cap E_{q(j)})\cup (D_{q(j)}\cap E_{p(i)})\not=\emptyset.$
\begin{claimproof}[Proof of claim]Just note that by the conditions (2) and (3), we have the following chain of equalities: \begin{align*}
   ( (D_{p(i)}\cap E_{q(j)})\cup (D_{q(j)}\cap E_{p(i)}))\cap m&=(D_{p(i)}\cap (E_{q(j)}\cap m))\cup ((D_{q(j)}\cap m)\cap E_{p(i)})\\&
    =(D_{p(i)}\cap (E_{p(j)}\cap m))\cup ((D_{p(j)}\cap m)\cap E_{p(i)})\\& 
    =(D_{p(i)}\cap E_{p(j)})\cup (D_{p(j)}\cap E_{p(i)})\not=\emptyset.&    \end{align*}
    This proves the claim.    
\end{claimproof}

\noindent
\underline{Claim 2}: Let $\mathcal{A}'\in [\mathcal{A}]^{\omega_1}$ and $i<n$. Then there is $\mathcal{X}\in [\mathcal{A}']^{\omega_1}$ so that for any $k\in \omega$, if  $\{p\in \mathcal{X}\,:\,k\in D_{p(i)}\}$  is non-empty, then it is uncountable. Analogously with the set $\{p\in \mathcal{X}\,:\,k\in E_{p(i)}\}$.
\begin{claimproof}[Proof of claim] Let $M$ be a countable elementary submodel of $H(\omega_2)$ such that $\mathcal{A}'\in M$. We put $\mathcal{X}= \mathcal{A}'\backslash M$. Note that if $k\in \omega$ is such that $k\in D_{q(i)}$ for some $q\in \mathcal{X}$, then the set $\{p\in \mathcal{A}'\,:\,k\in D_{p(i)}\}$ is an element of $M$ which is not contained in it. By elementarity, it follows that this set is uncountable. Therefore  $\{p\in \mathcal{X}\,:\,k\in D_{p(i)}\}=\{p\in \mathcal{A}'\,:\,k\in D_{p(i)}\}\backslash M$ is uncountable as well. The same argument holds for the set $\{p\in \mathcal{X}\,:\,k\in E_{p(i)}\}$. Hence, the proof is over.
\end{claimproof}
\noindent
\underline{Claim 3}: Let $i<n$ and $\mathcal{X}$, $\mathcal{Y}$ be uncountable disjoint subsets of $\mathcal{A}$ satisfying the conclusions of the Claim 2 when applied to $i$. Then there are $\mathcal{X}'\in [\mathcal{X}]^{\omega_1}$ and $\mathcal{Y}'\in [\mathcal{Y}]^{\omega_1}$ so that $$(D_{p(i)}\cap E_{q(i)})\cup(D_{q(i)}\cap E_{p(i)})\not=\emptyset$$
for all $p\in \mathcal{X}'$ and $q\in \mathcal{Y}'$.
\begin{claimproof}[Proof of claim.]Since $(\mathcal{D},\mathcal{E})$ is almost Luzin, we have that $(\langle D_{p(i)}\rangle_{p\in \mathcal{X}}, \langle E_{q(i)}\rangle_{q\in \mathcal{Y}})$ is a gap. In particular, there is $k\in \omega$, $p'\in \mathcal{X}$ and $q'\in \mathcal{Y}$ for which $k\in (D_{p'(i)}\cap E_{q'(i)})\cup(D_{q'(i)}\cap E_{p'(i)})$. Without loss of generality we may assume that $k\in D_{p'(i)}\cap E_{q'(i)}$. We define $$\mathcal{X}'=\{p\in \mathcal{X}\,:,k\in D_{p(i)}\},$$
$$\mathcal{Y}'=\{q\in \mathcal{Y}\,:\,k\in E_{q(i)}\}.$$
It is straightforward that $\mathcal{X}'$ and $\mathcal{Y}'$ are the sets that we are looking for.
\end{claimproof}
\noindent
By applying multiple times the claims  2 and 3, we may build two sequences $\mathcal{X}_{n-1}\subseteq\dots \subseteq \mathcal{X}_0\subseteq \mathcal{A}$ and $\mathcal{Y}_{n-1}\subseteq \dots \subseteq \mathcal{Y}_0\subseteq \mathcal{A}$ of uncountable sets for which $\mathcal{X}\cap \mathcal{Y}=\emptyset$ and such that $(D_{p(i)}\cap E_{q(i)})\cup(D_{q(i)}\cap E_{q(i)})\not=\emptyset$
for all $i<n$, $p\in \mathcal{X}_i$ and $q\in \mathcal{Y}_i$. Using Claim 1, it should be clear that if $p\in \mathcal{X}_{n-1}$ and $q\in \mathcal{Y}_{n-1}$, then $p\cup q$ is a condition of $\mathbb{L}(\mathcal{D},\mathcal{E})$. Thus, we are done.
\end{proof}
\end{proposition}
\begin{corollary}\label{corollarycccbiorthogonal}Let $(\mathcal{D},\mathcal{E})$ be an almost Luzin gap.  Then there is a $ccc$ forcing $\mathbb{Q}$ so that $\mathbb{Q}\Vdash\text{\say{ $(\mathcal{D},\mathcal{E})$ has a biorthogonal subgap.}}$
\end{corollary}
\begin{theorem}[Under PFA]\label{stronglyseppfa}Every $(\omega_1,\omega_1)$-gap is weakly levelwise-inseparable.
\begin{proof} Let $(\mathcal{L},\mathcal{R})=(L_\alpha,R_\alpha)_{\alpha\in\omega_1}$ be an $(\omega_1,\omega_1)$-gap. There is no loss of generality in assuming that $(\mathcal{L},\mathcal{R})$ is normal. Consider the forcing $\mathbb{P}=2^{<\omega_1}$ and let $G\subseteq \mathbb{P}$ be a $\mathbb{P}$-generic filter over $V$.
Then $V[G]$ models CH. Furthermore, $[\omega]^{\omega}\cap V=[\omega]^{\omega}\cap V[G]$. In this way, $(\mathcal{L},\mathcal{R})$ is still a gap in $V[G]$. According to the Proposition \ref{highlyadequateprop}, we can find a club $X\subseteq \omega_1$ in $V[G]$ which is highly adequate for $(\mathcal{L},\mathcal{R})$. That is, in $V[G]$, the gap $(L_{X(\alpha+1)}\backslash L_{X(\alpha)},R_{X(\alpha+1)}\backslash R_{X(\alpha)})_{\alpha\in \omega_1}$ is almost Luzin. Therefore, by virtue of the Corollary \ref{corollarycccbiorthogonal} there is a $ccc$ forcing $\mathbb{Q}$ which forces this gap to have a biorthogonal subgap. Let $H$ be a $\mathbb{Q}$-generic filter over $V[G]$. In $V[G][H]$, there is $Y\in [\omega_1]^{\omega_1}$ so that the pregap $$(\langle L_{X(\alpha+1)}\backslash L_{X(\alpha)}\rangle, \langle R_{X(\alpha+1)}\backslash R_{X(\alpha)}\rangle)_{\alpha\in Y}$$
is biorthogonal.\\
Returning to $V$, let $\dot{\mathbb{Q}}$ be a $\mathbb{P}$ name for $\mathbb{Q}$ in $V[G]$. Also, let $\dot{X}$ and $\dot{Y}$ be two $\mathbb{P}*\dot{\mathbb{Q}}$-names for $X$ and $Y$ respectively. Finally, let $p\in \mathbb{P}*\dot{\mathbb{Q}}$ which forces $\dot{X}$ and $\dot{Y}$ to have the properties discussed in the previous paragraph. Given $\xi\in \omega_1$, we define $$D_\xi:=\{q\leq p\:\,:\exists \xi<\alpha,\beta\in \omega_1\,(q\Vdash\text{\say{ $\alpha\in\dot{X}$, and $\beta\in \dot{Y}$)}}\,\},$$
$$E_\xi:=\{q\leq p\,:\,p\Vdash\text{\say{$\xi\not\in\dot{X}$}}\textit{ or there is }\alpha\in \omega_1\textit{ such that }p\Vdash\text{\say{$\dot{X}(\alpha)=\xi$}}\}.$$
It is straightforward that both $D_\xi$ and $E_\xi$ are dense below $p$ for any $\xi\in\omega_1$. Since $\mathbb{P}$ is $\sigma$-closed and $\dot{\mathbb{Q}}$ is forced to be $ccc$, then $\mathbb{P}*\dot{\mathbb{Q}}$ is proper. In this way, there is a filter $F$ which intersects  each $D_\xi$ and $E_\xi$. Let $$X'=\{\alpha\,:\,\exists q\in F\,(q\Vdash\text{\say{$\alpha\in \dot{X}$}})\},$$
$$Y'=\{\alpha\,:\,\exists q\in F\,(q\Vdash\text{\say{$\alpha\in \dot{Y}$}})\}.$$
Note that the statement \begin{center}
    \say{$(L_{X'(\alpha+1)}\backslash L_{X'(\alpha)}\cap R_{X'(\beta+1)}\backslash R_{X'(\beta)})\cup  (R_{X'(\alpha+1)}\backslash R_{X'(\alpha)}\cap L_{X'(\beta+1)}\backslash L_{X'(\beta)})=\emptyset$}
\end{center}
is absolute. Hence, it follows that that the pregap $$(L_{X'(\alpha+1)}\backslash L_{X'(\alpha)}, R_{X'(\alpha+1)}\backslash R_{X'(\alpha)})_{\alpha\in Y'}$$ is biorthogonal. In particular, it is a gap by Lemma \ref{biorthogonalgaplemma}. Form this, we  easily have that $(\mathcal{L},\mathcal{R})$ is weakly levelwise-inseparable.
\end{proof}
\end{theorem}

 \section{Basics on construction schemes}\label{sectionschemes}
We now recall what construction schemes are, and review some important facts about them. For a detailed introduction to this topic, see \cite{cruztesis}. A \emph{type} is a sequence $\tau=\langle m_k,n_{k+1},r_{k+1}\rangle_{k\in\omega}$ of triplets of natural numbers such that the following conditions hold for any $k\in\omega$:\\
\begin{enumerate}[label=$(\alph*)$,itemsep=0.5em]
\begin{minipage}{5cm}
\item $m_0=1,$
\item $n_{k+1}\geq 2,$
\end{minipage}
\begin{minipage}{5cm}
\item $m_k>r_{k+1}, $
\item $m_{k+1}=r_{k+1}+(m_k-r_{k+1})n_{k+1}.$
\end{minipage}\\
\end{enumerate}
 
We say that type is \textit{good} if for each $r\in\omega$, there are infinitely many $k$'s for which $r=r_k$. Throughout this text, we will only work with good types. We say that  partition of $\omega$, namely $\mathcal{P}$, is \emph{compatible with $\tau$} if $\tau$ is good when restricted to each member of $\mathcal{P}$. That is, for each $P\in \mathcal{P}$ and every $r\in \omega$ there are infinitely many $k\in P$ for which $r_k=r$.

 A construction scheme (of type $\tau$) over $X$ is a family $\mathcal{F}\subseteq \text{Fin}(\omega_1)$ which: Is cofinal in $(\text{Fin}(\omega_1),\subseteq )$, any member of $\mathcal{F}$ has cardinality $m_k$ for some $k\in\omega$, and furthermore, if we put $\mathcal{F}_k:=\{ F\in \mathcal{F}\,:\,|F|=m_k\}$, then the following two properties are satisfied for each $k\in\omega$:\vspace{0.5em}
\begin{enumerate}[label=(\roman*),itemsep=0.5em]
\item $\forall F,E\in \mathcal{F}_k\big(\; E\cap F\sqsubseteq E,F\;\big)$,
\item $\forall F\in \mathcal{F}_{k+1}\;\exists F_0,\dots,F_{n_{k+1}-1}\in \mathcal{F}_k$ such that $$F=\bigcup\limits_{i<n_{k+1}}F_i.$$
Moreover, $\langle F_i\rangle_{i<n_{k+1}}$ forms a $\Delta$-system with root $R(F)$ such that  $|R(F)|=r_{k+1}$ and $R(F)<F_0\backslash R(F)<\dots < F_{n_{k+1}-1}\backslash R(F).$
\end{enumerate}\vspace{0.5em}
Note that for a given $F\in \mathcal{F}_{k+1}$, each of the $F_i$'s mentioned above can be written as $F[r_{k+1}]\cup F[\,[a_i,a_i+1)\,]$ where $a_i=r_{k+1}+i\cdot(m_{k+1}-m_k)$. In particular, this is saying that the family $\langle F_i\rangle_{i<n_{k+1}}$ is unique, so  we call it the \emph{canonical decomposition} of $F$.

To each construction scheme $\mathcal{F}$, we associate with some natural functions of countable codomain. The first of such functions is the one derived from the cofinality condition in the definition of a scheme. Namely, given $\alpha,\beta\in \omega_1$ there is $F\in \mathcal{F}$ such that $\{\alpha,\beta\}\subseteq F$. Hence, we can define $\rho:\omega_1^2\longrightarrow \omega$ as:$$\rho(\alpha,\beta)=\min (k\in\omega\,:\,\exists F\in \mathcal{F}\,(\{\alpha,\beta\}\subseteq F )\,).$$
For each finite $A\subseteq \omega_1$, we also define $$\rho^A=\max(\rho[A^2])=\max(\rho(\alpha,\beta)\,:\,\alpha,\beta\in A).$$
It is not hard to see that $\rho^F=n$ for each $F\in \mathcal{F}_n$.
The most important feature of the function $\rho$ is that it is an ordinal metric\footnote{Ordinal metrics were introduced by the third author in \cite{partitioningpairs}. See \cite{Walksonordinals} for a full introduction.}. This means that it satisfies the properties stated in the following lemma.
\begin{lemma}Let $\mathcal{F}$ be a construction scheme. The following properties hold for any $\alpha,\beta,\gamma\in dom(\mathcal{F})$ and each $k\in\omega$:
\begin{enumerate}[label=$(om_{\arabic*})$]
\item $\rho(\alpha,\beta)=0$ if and only if $\alpha=\beta$.
\item $\rho(\alpha,\beta)=\rho(\beta,\alpha).$ 
\item If $\alpha\leq \min(\beta,\gamma)$, then $\rho(\alpha,\beta)\leq \max(\,\rho(\alpha,\gamma),\rho(\beta,\gamma)\,)$.
\item $\{\xi\leq\alpha\,:\,\rho(\alpha,\xi)\leq k\}$ is finite.
\end{enumerate}
\end{lemma}
Given $\alpha\in \omega_1$ and $k\in\omega$ we can define the \emph{$k$-closure} of $\alpha$ as $(\alpha)_k:=\{\xi\leq \alpha\,:\rho(\alpha,\xi)\leq k\}$ and $(\alpha)_k^-:=(\alpha)_k\backslash \{\alpha\}$. Note that property $(om_4)$ is saying that all the $k$-closures are finite. It is a useful fact that for any $k\in\omega$ (with $m_k\leq |dom(\mathcal{F})|$) there is at least one $F\in \mathcal{F}_k$ such that $\alpha\in F$. Even more, for any such $F$ we have the equalities: $$F\cap (\alpha+1)=(\alpha)_k$$
$$F\cap \alpha=(\alpha)^-_k.$$
Given $\alpha\in dom(\mathcal{F})$, the $k$-cardinality function $\lVert \alpha\rVert_{\_}:\omega\longrightarrow \omega$ is defined as:
$$\lVert \alpha\rVert_k=|(\alpha)^-_k|.$$
It is helpful to remember that if $F\in \mathcal{F}_k$ is such that $\alpha\in F$, then $F(\lVert \alpha \rVert_k)=\alpha$. Recall that if $f,g:\omega\longrightarrow\omega$ are distinct,  we can define $\Delta(f,g):=\min (\,k\in \omega\,:\,f(k)\not=g(k)\,)$. We put $\Delta(f,g):=\omega$ whenever $f=g$. In the case of  construction schemes, this leads to the definition of the function $\Delta:\omega_1^2\longrightarrow \omega+1$ which is defined as: $$\Delta(\alpha,\beta)=\Delta(\lVert \alpha\rVert_{\_},\lVert \beta\rVert_{\_}).$$
In the following lemma, we describe the most basic properties of this function as well as its relation with the $\rho$-function.
\begin{lemma}Let $\mathcal{F}$ be a construction scheme and let $\alpha,\beta,\delta\in \omega_1$. Then the following properties hold:
\begin{enumerate}[label=$(dp_{\arabic*})$, itemsep=0.5em]
    \item If $\alpha<\beta$, then $\lVert \alpha\rVert_k<\lVert \beta\rVert_k$ for each $k\geq \rho(\alpha,\beta)$. In particular, $\Delta(\alpha,\beta)\leq \rho(\alpha,\beta)$ whenever $\alpha\not=\beta$.
    \item If $\Delta(\alpha,\beta)<\Delta(\beta,\delta)$, then $\Delta(\alpha,\beta)=\Delta(\alpha,\delta)$.
\end{enumerate}
\end{lemma}
The last canonical function that we will present is the one which is related to the property $(ii)$ in the definition of a scheme. Namely, given $\alpha\in \omega_1$ we define the function $\Xi_\alpha:\omega\longrightarrow \{-1\}\cup \omega$ as follows;  If $1\leq k\in\omega$, then there is $F\in \mathcal{F}_k$ so that $\alpha\in F$. According to the property $(ii)$ in the definition of a construction scheme, we have that either $\alpha\in R(F)$ or there is a unique $i<n_k$ such that $\alpha\in F_i\backslash R(F)$. We then define $$\Xi_\alpha(k):=\begin{cases}-1&\text{ if }\alpha\in R(F)\\
i&\text{ if }\alpha\in F_i\backslash R(F)\end{cases}$$ It can be proved that this definition does not depend on the choice of $F.$ Now, if $k=0$ we simply define $\Xi_\alpha(k)$ as $0$.

The  two lemmas below relate the functions $\rho$, $\Xi$ and $\Delta$.
\begin{lemma}\label{lemmahdeltarhoinequalities}Let $\alpha,\beta\in \omega_1$ be distinct ordinals and $k<\Delta(\alpha,\beta)$. If $h:(\alpha)_k\longrightarrow (\beta)_k$ is the only increasing bijection and $\delta,\gamma\in (\alpha)_k$ are such that $\delta\leq \gamma$ and $h(\delta)\not=\delta$ then the following happens:
\begin{enumerate}[label=$(\alph*)$]
\item$h(\gamma)\not=\gamma,$
\item$\rho(\alpha,\beta)\geq \rho(\gamma,h(\gamma))\geq \rho(\delta,h(\delta))\geq\Delta(\delta,h(\delta)) \geq \Delta(\gamma,h(\gamma))\geq \Delta(\alpha,\beta).$
\end{enumerate}
\end{lemma}
\begin{lemma}\label{lemmaxi}Let $\mathcal{F}$ be a construction scheme, $\alpha<\beta\in dom(\mathcal{F})$ and $1\leq k\in \omega$. Then:\vspace{0.5em}
\begin{enumerate}[label=$(\alph*)$]
\item If $k<\Delta(\alpha,\beta)$, then  $\Xi_\alpha(k)=\Xi_\beta(k).$
\item If $k=\rho(\alpha,\beta)$, then $0\leq \Xi_\alpha(k)<\Xi_\beta(k).$
\item If $k>\rho(\alpha,\beta)$, then either $\Xi_\alpha(k)=-1$ or $\Xi_\alpha(k)=\Xi_\beta(k).$ 
\item If $k=\Delta(\alpha,\beta)$ then $0\leq \Xi_\alpha(k)\not=\Xi_\beta(k)\geq 0.$ 
\end{enumerate}
\end{lemma}

\begin{definition}[Captured systems] Let $1\leq n,m\in\omega$ and $\mathcal{D}=\langle D_i\rangle _{i<n}\subseteq [X]^{m}$ be a root-tail-tail $\Delta$-system with a root $R$ of cardinality $r$. Given $l\in \omega$, we will say that $\mathcal{D}$ is \emph{captured} at level $l$ if :\begin{itemize}
\item
For all $i<n$ and each $a<m$, $$\Xi_{D_i(a)}(l)=\begin{cases}-1&\text{if }a\leq r\\
i&\text{if }a>r
\end{cases}$$
 \item For all $i<j<n$ and $r\leq a<m$,
$$\rho(D_i(a),D_j(a))=l=\Delta(D_i(a),D_j(a)).$$
\end{itemize}
If moreover $n=n_l$, we say that $\mathcal{D}$ is fully captured.  Whenever $D\in \Fin(X)$, we say that $D$ is captured (resp. fully captured) in case $\{\{\alpha\}\,:\,\alpha\in D\}$ is captured (resp. fully captured).
\end{definition}

A scheme $\mathcal{F}$ is said to be
\emph{$\mathcal{P}$-$n$-capturing} (for some $n\in\omega$ and a partition $\mathcal{P}$  of $\omega$) if for any uncountable $\mathcal{S}\subseteq \text{Fin}(\omega_1)$ and each $P\in\mathcal{P}$, there are infinitely many $l\in P$ for which there is $\mathcal{D}\in [S]^n$ captured at level $l.$  $\mathcal{F}$ is $\mathcal{P}$ capturing provided that it is $\mathcal{P}$-$n$-capturing for every $n\geq 2$. The notion of a $\mathcal{P}$-fully capturing scheme has the expected meaning.  If $\mathcal{P}=\omega$, we omit it from the notation.

We will frequently use the following results.

\begin{lemma}\label{capturedfamiliestosetslemma}Let $\mathcal{D}$ be a finite subset of $\text{Fin}(X)$ and $l\in \omega$.  For any $A\in \mathcal{D}$, let $\alpha_A=\max(A)$. Suppose that there are $j\in \omega$, $a<m_j$ and $C\subseteq a+1$ so that the following conditions hold for any $A\in \mathcal{D}$:
\begin{enumerate}[label=$(\alph*)$]
    \item $\rho^A\leq j,$
    \item $\lVert \alpha_A\rVert_j=a$,
    \item $(\alpha_A)_j[C]=A.$ 
\end{enumerate}
If  $D=\{\alpha_A\,:\,A\in \mathcal{D}\}$ is captured at some level $l$, then $\mathcal{D}$ is also captured at level $l$.
\end{lemma}
\begin{lemma}\label{equivalencecapturing}Let $\mathcal{F}$ be a construction scheme and $\mathcal{P}$ be a partition of $\omega$ compatible with $\tau$. Then:
\begin{itemize}
\item For each $n\in\omega$, $\mathcal{F}$ is $n$-$\mathcal{P}$-capturing if and only if for each $S\in[\omega_1]^{\omega_1}$ and $P\in \mathcal{P}$ there is $D\in [S]^n$ which is captured at some level $l\in P$.
\item $\mathcal{F}$ is $\mathcal{P}$-fully capturing if and only if for each $S\in [\omega]^{\omega_1}$ and $P\in \mathcal{P}$ there is $D\in \text{Fin}(S)$ which is fully captured at some level $l\in P$.
\end{itemize}
\end{lemma}

We are now ready to state the Capturing axioms. The axiom FCA(part) was introduced in \cite{schemenonseparablestructures}. The other axioms were later studied in 
\cite{irredundantsetsoperator},
\cite{forcingandconstructionschemes}, \cite{banachspacescheme}, \cite{lopezschemethesis} and \cite{schemescruz}. All of them hold under $\Diamond$ or after adding at least $\omega_1$-Cohen reals (We do not know whether one Cohen real is enough).\\\\
{\bf Fully Capturing Axiom} {[\bf FCA]}: There is a fully capturing construction scheme over $\omega_1$ of every possible good type.\\\\
{ \bf Fully Capturing Axiom with Partitions [FCA(part)]}: There is a $\mathcal{P}$-fully capturing construction scheme over $\omega_1$ for every good type $\tau$  and each partition $\mathcal{P}$ compatible with $\tau$.\\\\
{\bf $n$-Capturing Axiom [CA$_n$]}: For any $n'\leq n$, there is an $n'$-capturing construction scheme over $\omega_1$ of every possible good type satisfying that $n'\leq n_k$ for each $k\in \omega\backslash1$.\\\\
{\bf $n$-Capturing Axiom with Partitions [CA$_n$(part)]}: For any $n'\leq n$, there is a $\mathcal{P}$-$n'$-capturing construction scheme over $\omega_1$ for every good type $\tau$ satisfying that $n'\leq n_k$ for each $k\in \omega\backslash 1$  and each partition $\mathcal{P}$ compatible with $\tau$.\\\\
{\bf Capturing Axiom [CA]}: CA$_n$ holds for each $n\in\omega$  and there is a capturing construction scheme over $\omega_1$ for every good type satisfying that the sequence $\langle n_{k+1}\rangle_{k\in\omega}$ is non-decreasing and unbounded.\\\\
{\bf Capturing Axiom with partitions [CA(part)]}: CA$_n$(part) holds for each $n\in\omega$ and there is a $
\mathcal{P}$-capturing construction scheme over $\omega_1$ for every good type $\tau$ satisfying that the sequence $\langle n_{k+1}\rangle_{k\in\omega}$ is non-decreasing and unbounded and each partition $\mathcal{P}$-compatible with $\tau$.\\

We now define the main concept of this work. Let us say that a forcing $\mathbb{P}$ \emph{$n$-preserves} a scheme $\mathcal{F}$ if $\mathbb{P}\Vdash\text{\say{$\mathcal{F}\textit{ is }n-capturing$}}$.
\begin{definition}[Parametrized Martin's numbers]\label{parametrizedmartinsdef}Let $\mathcal{F}$ be a construction scheme and $2\leq n\in\omega$. We define $\mathfrak{m}^n_\mathcal{F}$ as follows:
$$\mathfrak{m}^n_\mathcal{F}=\begin{cases}
\omega &\textit{if }\mathcal{F}\textit{ is not }n\textit{-capturing}\\
\min(\mathfrak{m}(\mathbb{P})\,:\,\mathbb{P}\textit{ is }ccc\textit{ and }\mathbb{P}\Vdash \text{\say{$\mathcal{F}\textit{ is }n-capturing$}})&\textit{if }\mathcal{F}\textit{ is }n\textit{-capturing}
\end{cases}$$
$\mathfrak{m}^2_\mathcal{F}$ is denoted simply as $\mathfrak{m}_\mathcal{F}.$ 
\end{definition}

\section{Not every gap is weakly levelwise-inseparable}\label{noteverygapiswli}

 In 1909,  Hausdorff gave a clever recursive construction of an $(\omega_1,\omega_1)$-gap without appealing to any extra axioms. This is quite surprising, as it is required to overcome $\mathfrak{c}$ obstacles in only $\omega_1$-many steps. Such gap satisfied the property below.

\begin{definition}[Hausdorff condition]Let $(\mathcal{L},\mathcal{R})=(L_\alpha,R_\alpha)_{\alpha\in\omega_1}$ be an $(\omega_1,\omega_1)$-pregap on $\omega$. We say that $(\mathcal{L},\mathcal{R})$ is  \emph{Hausdorff} if $\{\alpha<\beta\,:\,L_\beta\cap R_\alpha\subseteq k\}$ is finite for each $\beta\in \omega_1$ and $k\in\omega.$
\end{definition}
It is a standard excercise to prove that any $(\omega_1,\omega_1)$-pregap satisfying the Hausdorff condition is in fact a gap. Under the $P$-ideal Dichotomy (PID), if $(\mathcal{L},\mathcal{R})$ is an $(\omega_1,\omega_1)$-gap, then there are cofinal $\mathcal{L}'\subseteq \mathcal{L}$ and $\mathcal{R}'\subseteq \mathcal{R}$ so that $(\mathcal{L}',\mathcal{R}')$ is a Hausdorff gap. This was proved by Abraham and the third author in \cite{partitionpropertiesch}. The same result follows from $\mathfrak{m}>\omega_1.$

\begin{theorem}\label{hausdorffgapconstruction}Let $\mathcal{F}$ be a $2$-construction scheme. For each $\alpha\in \omega_1$, define $$L_\alpha=\{2k+\Xi_\alpha(k)\,:\,k\in \omega\backslash 1,\textit{ and }\Xi_\alpha(k)\geq 0\},$$
$$R_\alpha=\{2k+(1-\Xi_\alpha(k))\,:\,k\in\omega\backslash 1,\textit{ and }\Xi_\alpha(k)\geq 0\}.$$
Then, $(L_\alpha, R_\alpha)_{\alpha\in\omega_1}$ is a Hausdorff gap.
\begin{proof} Since $r_k=0$ for infinitely many $k's$, it should be clear that each $L_\alpha$ and $R_\alpha$ are infinite.  By virtue of the part (c) in Lemma \ref{lemmaxi}, we have that if $\alpha<\beta$ then  $$L_\alpha\backslash L_\beta\subseteq \{2k+\Xi_\alpha(k)\,:\,k\leq \rho(\alpha,\beta)\textit{ and }\Xi_\alpha(k)\not=\Xi_\beta(k)\,\},$$
$$R_\alpha\backslash R_\beta\subseteq \{2k+(1-\Xi_\alpha(k))\,:\,k\leq \rho(\alpha,\beta)\textit{ and }\Xi_\alpha(k)\not=\Xi_\beta(k)\,\}.$$ 
As the sets on the right are finite, this means  $\langle L_\alpha\rangle_{\alpha\in \omega_1}$ and $\langle R_\alpha\rangle_{\alpha\in \omega_1}$ are both $\omega_1$-towers. Furthermore by definition we have that $L_\alpha\cap R_\alpha=\emptyset$ for each $\alpha\in\omega_1$. This  implies that $(L_\alpha,R_\alpha)_{\alpha\in\omega_1}$ is a pregap.

The only thing left to prove is that the Hausdorff condition is satisfied. For this purpose take $\beta\in \omega_1$ and $k\in\omega$. We claim  $\{\alpha<\beta\,:\,L_\beta\cap R_\alpha\subseteq k\}\subseteq (\beta)_k$. For this, take an arbitrary $\alpha<\beta$ satisfying $\rho(\alpha,\beta)\geq k.$ By the part (b) of Lemma \ref{lemmaxi},  $\Xi_\alpha(\rho(\alpha,\beta))=0$ and $\Xi_\beta(\rho(\alpha,\beta))=1$. This means $2\rho(\alpha,\beta)+1\in L_\beta\cap R_\alpha$, so we are done.
\end{proof}
\end{theorem}
An interesting feature of the proof given above is that not only we prove that there is a gap using  a $2$-construction scheme, but such gap can be explicitly defined from it.

\begin{theorem}\label{firsthausdorffdonutseparable}The Hausdorff gap $(\mathcal{L},\mathcal{R})$ constructed in Theorem \ref{hausdorffgapconstruction} is levelwise-separable.
\begin{proof}Let $\alpha\in\omega_1$ and consider an arbitrary $k> \rho(\alpha,\alpha+1)$. Note that $\lVert \alpha\rVert_k +1=\lVert \alpha+1 \rVert_k$. By the point (c) of Lemma \ref{lemmaxi} we know that either $\Xi_{\alpha+1}(k)=\Xi_{\alpha}(k)$ or these two numbers are different and $\Xi_{\alpha}(k)=-1$. It follows that in the latter case $r_k-1=\lVert \alpha\rVert_k<\lVert \alpha+1\rVert_k=r_k<m_k$, which means that $\Xi_{\alpha+1}(k)=0$.  From this fact we conclude that \begin{align*}L_{\alpha+1}\backslash L_{\alpha}&=^*\{2k+\Xi_{\alpha+1}(k)\,:k>\rho(\alpha,\alpha+1)\textit{ and }\Xi_{\alpha}(k)\not=\Xi_{\alpha+1}(k)\}\\
&\subseteq \{ 2k\,: k\in\omega\textit{ and }\Xi_{\alpha+1}(k)=0\}\subseteq \{2k\,: k\in\omega\}.
\end{align*}
In the same way, we have that $R_{\alpha+1}\backslash R_\alpha\subseteq ^* \{2k+1\,:k\in\omega\}.$
Hence, the set of even numbers separates $(L_{\alpha+1}\backslash L_\alpha,R_{\alpha+1}\backslash R_\alpha)_{\alpha\in\omega_1}$.
\end{proof}
\end{theorem}
 \begin{lemma}\label{separabilityequivalence} Let $(\mathcal{L},\mathcal{R})=(L_\alpha,R_\alpha)_{\alpha\in\omega_1}$ be an $(\omega_1,\omega_1)$-gap. The following statements are equivalent:
 \begin{enumerate}[label=$(\arabic*)$]
 \item $(L_\alpha,R_\alpha)_{\alpha\in\omega_1}$ is not weakly levelwise-inseparable. 
 \item For any $S\in[\omega_1]^{\omega_1}$, the gap $(L_\alpha,R_\alpha)_{\alpha\in S}$ is levelwise-separable.
 \item For any club $S\subseteq \omega_1$, the gap $(L_\alpha,R_\alpha)_{\alpha\in S}$ is levelwise-separable.
 \end{enumerate}
 
\begin{proof}The only nontrivial part of the proof is to show that $(3)$ implies $(1)$. For this, let $X$ and $Y$ be uncountable subsets of $\omega_1$. Our goal is to prove that there is a $C\in [\omega]^{\omega}$ which separates $(L_{X(\alpha+1)}\backslash L_{X(\alpha)},R_{Y(\alpha+1)}\backslash R_{Y(\alpha)})_{\alpha\in\omega_1}$. Let $\langle M_\alpha\rangle_{\alpha\in\omega_1}$ be a continuous chain of elementary submodels of a largely enough $H(\lambda)$ so that $X,Y, (\mathcal{L},\mathcal{R})\in M_0$. Now define $S$ as $\{0\}\cup \{ M_\alpha\cap \omega_1\,:\,\alpha\in \omega_1\}$. Then $S$ is a club. Furthermore, by elementarity it is straightforward that for any $\alpha\in\omega_1$ there is $\beta \in S$ with:
\begin{multicols}{2}\begin{itemize}
\item $S(\beta)\leq X(\alpha)<X(\alpha+1)\leq S(\beta+1)$.
\item $S(\beta)\leq Y(\alpha)<Y(\alpha+1)\leq S(\beta+1).$
\end{itemize}
\end{multicols}
According to the hypotheses, there is a subset $C$  of $\omega$ which separates the pregap $(L_{S(\alpha+1)}\backslash L_{S(\alpha)}, R_{S(\alpha+1)}\backslash R_{S(\alpha)})_{\alpha\in\omega_1}$. Note that if $\alpha\in\omega_1$ and $\beta\in S$ are as previously stated then $L_{X(\alpha+1)}\backslash L_{X(\alpha)}\subseteq^* L_{S(\beta+1)}\backslash L_{S(\beta)}\subseteq^* C$ and $R_{Y(\alpha+1)}\backslash R_{Y(\alpha)}\subseteq^* R_{Y(\beta+1)}\backslash R_{Y(\beta)}\subseteq^* \omega\backslash C$. Therefore, $C$ separates $$(L_{X(\alpha+1)}\backslash L_{X(\alpha)},R_{Y(\alpha+1)}\backslash R_{Y(\alpha)})_{\alpha\in\omega_1}.$$ This finishes the proof.
\end{proof}
\end{lemma}
\begin{corollary}\label{corocccpreservingdonut}Let $(\mathcal{L},\mathcal{R})=(L_\alpha,R_\alpha)_{\alpha\in\omega_1}$ be an $(\omega_1,\omega_1)$-gap which is not weakly levelwise-inseparable. If $\mathbb{P}$ is a $ccc$ forcing, then $$\mathbb{P}\Vdash\text{\say{ $(\mathcal{L},\mathcal{R})$ is not weakly levelwise-inseparable}}.$$
\begin{proof}Suppose towards a contradiction that there is $G$ a $\mathbb{P}$-generic filter over $V$ such that, in $V[G]$,  $(\mathcal{L},\mathcal{R})$ is weakly levelwise-inseparable.  According to the point (3) of Lemma \ref{separabilityequivalence}, there is a club $S\in V[G]$ such that the gap $(L_\alpha,R_\alpha)_{\alpha\in S}$ is levelwise-inseparable. Since $\mathbb{P}$ is a $ccc$ forcing, we can find a club $C\in V$ for which $C\subseteq S$. In particular, $C$ is levelwise-inseparable and belongs to $V$. This contradiction ends the proof. 
\end{proof}
    
\end{corollary}
Suppose that $(\mathcal{L},\mathcal{R})$ is a pregap over $\omega$ and $C$ is an infinite subset separating it. Define $s:\mathcal{L}\cup\mathcal{R}\longrightarrow \omega$ as: $$s_C(X)=\begin{cases}
\min(n\,: X\backslash n\subseteq C)&\textit{ if }X\in \mathcal{L}\\
\min(n\,:X\backslash n \cap C=\emptyset)&\textit{ if }X\in \mathcal{R}
\end{cases}$$
Then $L\cap R\subseteq \max(s_C(L),s_C(R))$ for any $L\in \mathcal{L}$ and $R\in \mathcal{R}.$ This motivates the following definition.
\begin{definition}[Separating functions]Let $(\mathcal{L},\mathcal{R})$ be a pregap. We say that $s:\mathcal{L}\cup\mathcal{R}\longrightarrow \omega$ is \textit{separating} if $L\cap R\subseteq  \max(s(L),s(R))$ for any $L\in \mathcal{L}$ and $R\in \mathcal{R}$. 
\end{definition}
\begin{lemma}Let $(\mathcal{L},\mathcal{R})$ be a pregap. If there is a separating $s:\mathcal{L}\cup\mathcal{R}\longrightarrow \omega$ then there is $C\in [\omega]^{\omega}$ which separates $(\mathcal{L},\mathcal{R})$.
\begin{proof}Let $s$ be as in the hypotheses. We define $C$ as $$\{n\in\omega\,:\, \exists L\in \mathcal{L}\,(n\in L\textit{ and }s(L)<n)\}=\bigcup\limits_{L\in \mathcal{L}}L\backslash (s(L)+1).$$
Note that $L\subseteq^* C$ for any $L\in \mathcal{L}$. For $R\in \mathcal{R}$ we claim that $R\backslash (s(R)+1)\cap C=\emptyset$ which in particular implies that $R\cap C=^* \emptyset$. Suppose towards a contradiction that this is not the case and let $n$ be an element in the intersection of both sets. On one hand, since $n\in C$ there is $L\in \mathcal{L}$ with $n\in L$ and $s(L)<n$. In particular, $n\in L\cap R$.  On the other hand, since $s$ is separating then $L\cap R\subseteq \max(s(R),s(L))\subseteq n$. This is a contradiction. Therefore the claim is true, which means that $C$ separates $(\mathcal{L},\mathcal{R})$. 
\end{proof}
\end{lemma}

\begin{definition}\label{separabilityforcingdonuts}Let $(\mathcal{D},\mathcal{E})=(D_\alpha,E_\alpha)_{\alpha\in\omega_1}$ be a  normal pregap. We define the forcing $\mathbb{P}(\mathcal{D},\mathcal{E})=\mathbb{P}(D_\alpha,E_\alpha)_{\alpha\in\omega_1}$ as the set of all functions $p;\omega_1\longrightarrow \omega$ with finite domain and such that $D_\alpha\cap E_\beta\subseteq  \max(p(\alpha),p(\beta))$ for all $\alpha,\beta\in dom(p)$. The order is given by $$p\leq q\text{ if and only if }q\subseteq p.$$
    
\end{definition}

\begin{rem}\label{remdenseseparating}If $(\mathcal{D},\mathcal{E})=(D_\alpha,E_\alpha)_{\alpha\in \omega_1}$ is as in the previous definition and $\beta\in \omega_1$ then $\mathcal{M}_\beta=\{p\in \mathbb{P}(\mathcal{D},\mathcal{E})\,:\,\beta\in dom(p)\,\}$ is dense in 
$\mathbb{P}(\mathcal{E},\mathcal{D})$. If $\mathcal{G}$ is a filter intersecting each $\mathcal{M}_\beta$ then the function $s:\mathcal{D}\cup \mathcal{E}\longrightarrow\omega_1$ given by: $$s(D_\alpha)=s(E_\alpha)=\bigcup\mathcal{G}(\alpha)$$
is well-defined and separating.
\end{rem}
The following proposition generalizes a well-known result of Kunen (we reiterate that our pregaps do not need to be linearly ordered).

\begin{proposition}\label{separabilitygapsprop} Let $(E_\alpha,D_\alpha)_{\alpha\in\omega_1}$ be a normal pregap. The following statements are equivalent:
\begin{enumerate}[label=$(\alph*)$]
\item $\mathbb{P}=\mathbb{P}(E_\alpha,D_\alpha)_{\alpha\in\omega_1}$ is $ccc$.
\item There is $\mathbb{Q}$ $ccc$ with $\mathbb{Q}\Vdash \text{\say{ $(E_\alpha,D_\alpha)_{\alpha\in\omega_1}\textit{ can be separated}$}}$.
\item There is $W$ a transitive model of ZFC extending $V$ with $\omega_1^W=\omega_1^V$ where $$W\models (E_\alpha,D_\alpha)_{\alpha\in\omega_1}\textit{ can be separated.}$$
\end{enumerate}
\begin{proof}Trivially $(a)$ implies $(b)$ and $(b)$ implies $(c)$. In order to prove that $(c)$ implies $(a)$ take an arbitrary uncountable subset $\mathbb{P}$ in $V$, say $\mathcal{A}$. For  any $p\in \mathcal{A}$ let $d_p=dom(p)$. Without loss of generality we can suppose that $\mathcal{A}$ is a $\Delta$-system with root $R$, there is $n\in\omega$ so that $|p|=n$ for every $p\in\mathcal{A}$ and $p(d_p(i))=q(d_q(i))$ for all $p,q\in \mathcal{A}$ and each $i<n$.

Let $W$ be as in the hypotheses of $(c)$ and let $C\in W$ which separates $(E_\alpha,D_\alpha)_{\alpha\in\omega_1}$. Since $\omega_1^V=\omega_1^W$ and $V$ models that $\mathcal{A}$ is uncountable, we can find an uncountable $\mathcal{A}'\subseteq \mathcal{A}$ in $W$, so that $E_{d_p(i)}\backslash C=E_{d_q(i)}\backslash C$ and $D_{d_p(i)}\cap C= D_{d_q(i)}\cap C$ for all $p,q\in \mathcal{A}'$ and $i<n$. Fix two distinct $p,q\in\mathcal{A}'$. We claim that $p\cup q\in \mathbb{P}$. Indeed, let $i,j<n$. Then \begin{align*}E_{d_p(i)}\cap D_{d_q(j)}&=((E_{d_p(i)}\backslash C)\cap D_{d_q(j)})\cup (E_{d_p(i)}\cap (D_{d_q(j)}\cap C))\\
&=((E_{d_q(i)}\backslash C)\cap D_{d_q(j)})\cup (E_{d_p(i)}\cap (D_{d_p(j)}\cap C))\\
&\subseteq \max(q(d_q(i)),q(d_q(j))\cup \max(p(d_p(i)),p(d_p(j)))\\
&= \max(p(d_p(i)),p(d_q(j))).
\end{align*}
This finishes the proof.
\end{proof}

\end{proposition}
Consider the gap $(L_\alpha, R_\alpha)_{\alpha\in\omega}$ constructed in Theorem \ref{hausdorffgapconstruction} and for any $k\in\omega$ let $N_k=2(k+1)=\{0,\dots,2k+1\}$. The following properties follow directly from the proof of such theorem and the definition of the $\Delta$-function.
\begin{proposition}\label{propositioninterhausdorf}Let $\alpha<\beta\in \omega_1$. Then:
\begin{enumerate}[label=$(\arabic*)$]
\item[$(0)$]$L_\alpha\cap R_\alpha=\emptyset.$
\item For each $k\in\omega$, both $\{2k,2k+1\}\cap L_\alpha$ and $\{2k,2k+1\}\cap R_\alpha$ have at most one point.  
\item All of the sets $L_\alpha\backslash L_\beta$,  $R_\alpha\backslash R_\beta$,  $L_\alpha\cap R_\beta$  and $L_\beta\cap R_\alpha$ are subsets of $N_{\rho(\alpha,\beta)}$.
\item $L_\alpha \cap N_{\Delta(\alpha,\beta)-1}=L_\beta\cap N_{\Delta(\alpha,\beta)-1}$ and  $R_\alpha \cap N_{\Delta(\alpha,\beta)-1}=R_\beta\cap N_{\Delta(\alpha,\beta)-1}$.

\end{enumerate}
\end{proposition}

In the next theorem, we will show that the existence of $(\omega_1,\omega_1)$-gaps which are not weakly levelwise-inseparable is consistent with ZFC.  The proof we provide can be simplified by the use of Lemma \ref{lemmaequivalencefunctioncapturingpreserving}. We decided to keep it this way in order to motivate such lemma. 
\begin{theorem}\label{mfstronglydonuttheorem} Let $\mathcal{F}$ be a $2$-capturing $2$-construction scheme for which $\mathfrak{m}_\mathcal{F}>\omega_1$. The Hausdorff gap constructed from $\mathcal{F}$ in Theorem \ref{hausdorffgapconstruction}  is not weakly levelwise-inseparable.
\begin{proof}Let $S\in [\omega_1]^{\omega_1}$. We will show that the gap $(L_{S(\alpha)},R_{S(\alpha)})_{\alpha\in\omega_1}$ is levelwise-separable. In other words, we will show that pregap  $(D_\alpha, E_\alpha)_{\alpha\in \omega_1}$ can be separated where $D_\alpha=L_{S(\alpha+1)}\backslash L_{S(\alpha)} $ and  $ E_\alpha=R_{S(\alpha+1)}\backslash R_{S(\alpha)}$ for each $\alpha\in \omega_1$. This is enough due to the Lemma \ref{separabilityequivalence}. In view of Remark \ref{remdenseseparating} and since we are assuming that $\mathfrak{m}_\mathcal{F}>\omega_1$, it is sufficient to prove that for $\mathbb{P}=\mathbb{P}(D_\alpha,E_\alpha)_{\alpha\in \omega_1}$ we have:
\begin{enumerate}[label=$(\arabic*)$]
\item $\mathbb{P}$ is a $ccc$ forcing,
\item $\mathbb{P}\Vdash\text{\say{ $ \mathcal{F}\textit{ is }2\textit{-capturing}$ }}.$
\end{enumerate}

\begin{claimproof}[Proof of $(1)$] For this consider an arbitrary $\mathcal{A}\in[\mathbb{P}]^{\omega_1}$. Given $p\in \mathcal{A}$ define $Z_p=\{ S(\alpha),S(\alpha+1)\,:\,\alpha\in dom(p)\}$. By refining $\mathcal{A}$ we can assume without loss of generality that for any $p,q\in \mathcal{A}$ the following conditions hold:
\begin{enumerate}[label=$(\alph*)$]
    \item $|p|=|q|$,
    \item If $f:dom(p)\longrightarrow dom(q)$ is the increasing bijection then $p(\alpha)=q(f(\alpha))$ for each $\alpha\in dom(p).$ In particular, there is $k\in\omega$ with $k>\max(im(p))$ for each $p\in \mathcal{A}$.
\end{enumerate}
Furthermore, we can suppose that the set $\{dom(p)\,:\,p\in \mathcal{A}\}$ forms a root-tail-tail $\Delta$-system with root $R$ satisfying the following properties for any two distinct $p,q\in \mathcal{A}$:
\begin{enumerate}
    \item[(c)] If $R\not=\emptyset$, then $\max(R)+1<\min(dom(p)\backslash R)$.
    \item [(d)] $\max(dom(p))+1\not\in dom(q)$.
\end{enumerate}
 As the set $\{Z_p\,:\,p\in \mathcal{A}\}$ is uncountable and $\mathcal{F}$ is assumed to be $2$-capturing then there are distinct $p,q\in \mathcal{A}$ for the set $\{Z_p,Z_q\}$ is captured at some level $l>k$. We affirm that $p$ and $q$ are compatible. This will follow from the next claim.\\\\
 \underline{Claim 1:} $p\cup q$ is a condition of $\mathbb{P}$.
 
 \begin{claimproof}[Proof of claim] First note that if $F\in \mathcal{F}_l$ is such that $Z_p\cup Z_q\subseteq F$ then $Z_p\cap Z_q=Z_p\cap R(F)=Z_q\cap R(F)$. 

Now, let us consider $f:dom(p)\longrightarrow dom(q)$  and $h:F_0\longrightarrow F_1$ the increasing bijections. In order to prove that $p\cup q$ is a condition it is enough to take $\alpha\in dom(p)$ and $\beta\in dom(q)$ and show that both $E_\alpha\cap D_\beta$ and $D_\alpha\cap E_\beta$ are contained in $\max(p(\alpha),q(\beta))$. If either $\alpha$ or $\beta$ belong to the intersection of $dom(p)\cap dom(q)$ there is nothing to do. So we can assume that $\alpha\in dom(p)\backslash dom(q)$ and $\beta\in dom(q)\backslash dom(p)$. 
Then $\{S(\alpha),S(\alpha+1)\} \subseteq Z_p\backslash Z_q\subseteq F_0\backslash R(F)$ and $\{S(\beta),S(\beta+1)\}\subseteq Z_q
\backslash Z_p\subseteq F_1\backslash R(F)$ due to the points (c) and (d).  Thus $\rho(S(\alpha+1),S(\beta+1))=l$ and consequently $D_\alpha\cap E_\beta\subseteq L_{S(\alpha+1)}\cap R_{S(\beta+1)}\subseteq N_l$ due to the point (2) of Proposition \ref{propositioninterhausdorf}.  But $\Xi_{S(\alpha)}(l)=\Xi_{S(\alpha+1)}(l)=0$. 
Therefore $2l\in L_{S(\alpha)}\cap L_{S(\alpha+1)}$. In particular $2l\in L_{S(\alpha+1)}$,  so $2l+1\not\in L_{S(\alpha+1)}$ by the point (1) of Proposition \ref{propositioninterhausdorf}. In this way  $\{2l,2l+1\}\cap D_\alpha=\emptyset$. Thus, $$D_\alpha\cap E_\beta\subseteq N_{l-1}.$$  
The next thing to note is that $h(S(\alpha))=S(f(\alpha))$ and $h(S(\alpha+1))=S(f(\alpha)+1).$ This means that $\Delta(S(\alpha),S(f(\alpha)))=l=\Delta(S(\alpha+1),S(f(\alpha)+1))$. In virtue of the point (3) of proposition \ref{propositioninterhausdorf} we have $$L_{S(\alpha)}\cap N_{l-1}=L_{S(f(\alpha))}\cap N_{l-1},$$
$$L_{S(\alpha+1)}\cap N_{l-1}=L_{S(f(\alpha)+1)}\cap N_{l-1}.$$
Hence, $D_\alpha\cap N_{l-1}=D_{f(\alpha)}\cap N_{l-1}$. From all the equations we have so far we deduce that $D_\alpha\cap E_\beta=D_{f(\alpha)}\cap E_\beta\subseteq \max(q(f(\alpha)),q(\beta))=\max(p(\alpha),q(\beta)).$ In a completely similar way we can show that $E_\alpha\cap D_\beta\subseteq \max(p(\alpha),q(\beta))$. We conclude that $p\cup q\in \mathbb{P}$ so we are done.
\end{claimproof}
\end{claimproof}

\begin{claimproof}[Proof of $(2)$]For this we will show that the equivalence of $2$-capturing stated in Lemma \ref{equivalencecapturing} is forced by $\mathbb{P}$. Let $\dot{X}$ be a name so that $\mathbbm{1}_\mathbb{P}\Vdash \text{\say{ $\dot{X}\in [\omega_1]^{\omega_1}$}}$ and let $q\in \mathbb{P}$. We shall find $p\leq q$ so that $$p\Vdash\text{\say{ $\exists M\in[\dot{X}]^2\,(\, M\textit{ is captured })$ }} .$$ As $\dot{X}$ is forced 
to be  uncountable we can find for each $\alpha\in \omega_1$ a condition $p_\alpha\leq q$ and $\xi_\alpha>\alpha$ with $p_\alpha\Vdash \text{\say{ $\xi_\alpha\in \dot{X}$ }}.$ We can suppose without loss of generality that the family $\mathcal{A}=\{p_\alpha\,:\,\alpha\in \omega_1\}$ satisfies the conditions (a), (b), (c) and (d) of the previous paragraphs. We may also assume that for any distinct $\alpha,\beta\in \omega_1:$
\begin{enumerate}
\item[$(d)$]$\xi_\alpha\not=\xi_\beta$ and $|Z_{p_\alpha}\cup \{\xi_\alpha\}|=|Z_{p_\beta}\cup \{\xi_\beta\}|$,
\item[$(e)$] If $g: Z_{p_\alpha}\cup \{\xi_\alpha\}\longrightarrow Z_{p_\beta}\cup \{\xi_\beta\}$ is the increasing bijection then $g[Z_{p_\alpha}]=Z_{p_\beta}$ and $g(\xi_\alpha)=\xi_\beta.$
\end{enumerate}
Since $\mathcal{F}$ is $2$-capturing there are $\alpha<\beta\in \omega_1$ for which $\{Z_{p_\alpha}\cup \{\xi_\alpha\},Z_{p_\beta}\cup \{\xi_\beta\}\}$ is captured at some level $l$ greater than $k$. In virtue of the conditions (d) and (e) it is easy to see that $\{Z_{p_\alpha},Z_{p_\beta}\}$ and $\{\xi_\alpha,\xi_\beta\}$ are also captured at level $l$. As in the first part of the proof, this implies that $p=p_\alpha\cup p_\beta$ is a condition of $\mathbb{P}$. To finish, just note that $p\leq q$ and  $p\Vdash\text{\say{ $\{\xi_\alpha,\xi_\beta\}\in [\dot{X}]^2\textit{ and }\{\xi_\alpha,\xi_\beta\}\textit{ is captured }$}}.$ 
\end{claimproof}
\end{proof}
\end{theorem}
The inequality $\mathfrak{m}_\mathcal{F}>\omega_1$ is consistent as we shall see in the following section. Therefore, 
by combining the results above and Theorem \ref{mfstronglydonuttheorem}, we conclude the following:
\begin{theorem}The statement \say{Every gap is weakly levelwise-inseparable} is independent from $ZFC.$ 
\end{theorem}
So far we have seen that both CH and PFA imply that every gap is weakly levelwise-inseparable. On the other hand, the existence of a $2$-capturing construction scheme $\mathcal{F}$ for which $\mathfrak{m}_\mathcal{F}>\omega_1$ implies the existence of a non WLI gap. We will end this subsection by proving that MA is not sufficient to decide this problem. 
\begin{theorem}\label{independentmastronglysep}MA is independent from the statement \say{Every gap is weakly levelwise-inseparable}.
\begin{proof}Since PFA implies that all gaps are WLI, then MA is consistent with that same statement.
In order to show the consistency of the negation, we start with a model $V$ in which there is a non WLI gap. We now consider a $ccc$ forcing $\mathbb{P}$ which forces  MA and let $G$ be a $\mathbb{P}$-generic filter over $V$. Then $V[G]$ is a model MA. Furthermore, by means of the Corollary \ref{corocccpreservingdonut}, $V[G]$ also models the existence of a non WLI gap. 
\end{proof}
\end{theorem}

\section{Fragments of Martin's axiom}\label{fragmentmartin}

Let us recall some of the previous work done by  Kaladjzievski and  Lopez \cite{forcingandconstructionschemes} regarding the preservation of capturing properties.

\begin{definition}[$n$-preserving property]Let $\mathbb{P}$ be a forcing notion and which preserves $\omega_1$ and $\mathcal{F}$ be a construction scheme. Given $2\leq n\in \omega$, we say that $\mathbb{P}$ \textit{$n$-preserves $\mathcal{F}$} if  $$\mathbb{P}\Vdash\text{\say{ $\mathcal{F}$ is $n$-capturing }}.$$
\end{definition}
Note that if $\mathbb{P}$ and $\mathcal{F}$ are as in the previous definition, then $\mathcal{F}$ is $n$-capturing.
 Usually, proving that a $ccc$ forcing $\mathbb{P}$ $n$-preserves a scheme $\mathcal{F}$, is similar to proving that $\mathbb{P}$ is $ccc$.  The following lemma gives us a useful way of handling this situation. 

\begin{lemma}\label{lemmaequivalencefunctioncapturingpreserving}Let $\mathcal{F}$ be an $n$-capturing construction scheme and $\mathbb{P}$ be a forcing. The two following statements are equivalent:
\begin{enumerate}[label=$(\arabic*)$]
    \item $\mathbb{P}$ is $ccc$ and $n$-preserves $\mathcal{F}.$
    \item For any $\mathcal{A}\in[\mathbb{P}]^{\omega_1}$ and each injective function $\nu:\mathcal{A}\longrightarrow \omega_1$, there is $\{p_0,\dots, p_{n-1}\}\in [\mathcal{A}]^n$ for which:
    \begin{itemize}
        \item $\{p_0,\dots, p_{n-1}\}$ is centered. That is, there is $p\in \mathbb{P}$ so that $p\leq p_i$ for any $i<n$.
        \item $\{\nu(p_0),\dots,\nu( p_{n-1})\}$ is captured.
    \end{itemize}
\end{enumerate}
\begin{proof}\begin{claimproof}[Proof of $(1)\Rightarrow(2)$]Let $\mathcal{A}$ be an uncountable subset of $\mathbb{P}$ and consider $\nu:\mathcal{A}\longrightarrow \omega_1$ an injective function. As $\mathbb{P}$ is $ccc$, there is $G$ a $\mathbb{P}$-generic filter over $V$ for which $G\cap \mathcal{A}$ is uncountable. In $V[G]$, let $X=\nu[\mathcal{A}\cap G]$. Note that $X$ is uncountable since $\nu$ is injective. As $\mathbb{P}$ $n$-preserves $\mathcal{F}$, there is $\{\alpha_0,\dots,\alpha_{n-1}\}\in [X]^n$ which is captured. Given $i<n$, let $p_i\in \mathcal{A}\cap G$ be such that $\nu(p_i)=\alpha_i$. Then $\{p_0,\dots,p_{n-1}\}$ is centered because it is included in $G$, and $\{\nu(p_0),\dots,\nu(p_{n-1})\}$ is captured. 
\end{claimproof}
\begin{claimproof}[Proof of $(1)\Leftarrow(2)$] In order to show that $\mathbb{P}$ is $ccc$ let $\mathcal{A}\in [\mathbb{P}]^{\omega_1}$. Let us consider $\nu:\mathcal{A}\longrightarrow \omega_1$ an aribtrary inyective function. According to the hypotheses, there is $\{p_0,\dots p_{n-1}\}\in [\mathcal{A}]^n$ which is centered. Particularly, $p_0$ and $p_1$ are two distinct compatible elements of $\mathcal{A}$. Thus, $\mathcal{A}$ is not an antichain.\\
Now we will show that $\mathbb{P}$ $n$-preserves $\mathcal{F}$. For this purpose, let $G$ be a $\mathbb{P}$-generic filter over $V$ and  $X\in V[G]\cap [\omega_1]^{\omega_1}$. Let $\dot{X}$ be a name for $X$ which is forced to be an uncountable subset of $\omega_1$ by $\mathbbm{1}_\mathbb{P}$. To finish, it is enough to show that the set of all $p\in \mathbb{P}$ for which there is $D\in [\omega_1]^{n}$ so that $D$ is captured and $p\Vdash\text{\say{ $D\subseteq \dot{X}$}}$, is dense in $\mathbb{P}$. For this purpose, let $q\in \mathbb{P}$. If there is $p\leq q$ for which the set $Y_p=\{\alpha\in \omega_1\,:\,p\Vdash\text{\say{$\alpha\in \dot{X}$} } \}$
is uncountable, we are done.  Therefore, we may assume that $Y_p$ is at most countable for each $p\leq q$. From this fact, it is easy to see that there is $\mathcal{A}\in[\mathbb{P}]^{\omega_1}$ and an injective function $\nu:\mathcal{A}\longrightarrow \omega_1$ so that $p\leq q$ and $p\Vdash \text{\say{ $\nu(p)\in \dot{X}$}}$ for any $p\in \mathcal{A}$. According to the hypotheses, there is $\{p_0,\dots,p_{n-1}\}\in [\mathcal{A}]^n$ which is centered and for which $D=\nu[\{p_0,\dots,p_{n-1}\}]$ is captured. Let $p$ be such that $p\leq p_i$ for each $i<n$. Then $p\leq q$ and $p\Vdash\text{\say{$D\subseteq \dot{X}$}}$. 
\end{claimproof}
    
\end{proof}
\end{lemma}
The lemma below was proved in \cite{forcingandconstructionschemes}.
\begin{lemma}\label{nknastercapturinglemma} Let $n\in \omega$, $\mathcal{P}$ a partition of $\omega$ and suppose that $\mathcal{F}$ is a $\mathcal{P}$-$n$-capturing construction scheme. If $\mathbb{P}$ is an $n$-Knaster forcing then $\mathbb{P}$ $\mathcal{P}$-$n$-preserves $\mathcal{F}$.
\end{lemma}
In \cite{forcingandconstructionschemes}, Kalajdzievski and Lopez used the  lemma above to show that $\mathfrak{m}_{K_n}$ is consistent with $CA_n$. 
    \begin{theorem}Let $2\leq n\in \omega$ and $V$ be a model of $CH$. For any cardinal $\kappa>\omega_1$ of uncountable cofinality, there is an $n$-Knaster forcing $\mathbb{K}$ so that $$\mathbb{K}\Vdash\text{\say{ $\mathfrak{c}=\kappa=\mathfrak{m}_{K_n}+CA_n$ }}.$$

\end{theorem}
In that same paper, they showed that $\mathfrak{m}_{K_n}>\omega_1$ implies that there are no $n+1$-capturing construction schemes. For that purpose, they used  the next property.\\\\
{\bf The Property $(\star)_n$:} For any uncountable $\Gamma\subseteq \omega^\omega$ there is $\Gamma_0\in[\Gamma]^{\omega_1}$ such that there are no $g_0,\dots,g_n\in \Gamma_0$ and $k\in \omega$ with $g_0|_k=\dots=g_n|_k$ and $g_i(k)\not=g_j(k)$ for any $i<j\leq n$.\\\\
The property $(\star)_n$ was considered by the third author in \cite{remarkscellularityproducts}. The next theorem follows from Lemma 6 of such paper.
\begin{theorem}Let $2\leq n\in \omega$. Then $\mathfrak{m}_{K_n}>\omega_1$ implies $(\star)_n$.
\end{theorem}
In contrast, $n+1$-capturing construction schemes imply the failure of $(\star)_n$ (see Theorem 2.4 in \cite{forcingandconstructionschemes}).
\begin{theorem}\label{starfailscapturing}Let $1\leq n\in \omega$. If there is an $n+1$-capturing construction scheme,  then $(\star)_n$ fails.
\begin{proof}Let $\mathcal{F}$ be an $n+1$-capturing construction scheme. For any $\alpha\in \omega_1$ let $f_\alpha:\omega\longrightarrow \omega$ be defined as:
$$f_\alpha(k)=\lVert \alpha\rVert_k.$$
Now, let $\Gamma=\{\,f_\alpha\,:\,\alpha\in \omega_1\}$. Since  $\mathcal{F}$ is $n+1$-capturing, for any $S\in [\omega_1]^{\omega_1}$ there is $D\in [S]^{n+1}$ which is captured at some level $l.$ Then $g_{D(0)}|_l=\dots=g_{D(n)}|_l$ and $g_{D(i)}(l)<g_{D(j)}(l)$ for any $i<j\leq n$. In this way, $\Gamma$ testifies the failure of $(*)_n$.
\end{proof}
\end{theorem}

\begin{rem} Suppose that $\mathcal{F}$ is an $n$-capturing construction scheme. The following diagram  represents the basic relations between the cardinals $\mathfrak{m}^i_\mathcal{F}$ and $\mathfrak{m}_{K_i}$ for each $i\leq n$. Note that the non-trivial relations hold due to Lemma \ref{nknastercapturinglemma}.
\begin{center}
\begin{tikzcd}
 &\mathfrak{m}_{K_2} \arrow[d, rightarrow] &\mathfrak{m}_{K_3} \arrow[l, rightarrow] \arrow[d, rightarrow]&\dots \arrow[l, rightarrow] \arrow[d, rightarrow] & \mathfrak{m}_{K_n}\arrow[l, rightarrow]\arrow[d, rightarrow] & \mathfrak{c} \arrow[l, rightarrow]\\
 \omega_1 & \mathfrak{m}^2_\mathcal{F} \arrow[l,rightarrow] & \mathfrak{m}^3_\mathcal{F} \arrow[l, rightarrow] &\dots\arrow[l, rightarrow] & \mathfrak{m}^n_\mathcal{F} \arrow[l, rightarrow] &
\end{tikzcd}
\end{center} 
\end{rem}
Having the $n$-Knaster property may be optimal in terms $n$-preserving a construction scheme.
\begin{proposition}\label{propequivalenceknasterpreservescheme}Let $2\leq n\in \omega$ and $\mathcal{F}$ be an $n$-capturing construction scheme for which $\mathfrak{m}_\mathcal{F}^n>\omega_1$. Given a $ccc$ forcing $\mathbb{P}$ and $m\geq n$, the following statements are equivalent:
\begin{enumerate}[label=$(\alph*)$]
\item $\mathbb{P}$ is $m$-Knaster.
    \item $\mathbb{P}$ is $n$-Knaster.
    \item $\mathbb{P}$ $n$-preserves $\mathcal{F}$.
    \item $\mathbb{P}$ has precaliber $\omega_1$.
\end{enumerate}
    \begin{proof} $(a)$ implies $(b)$ and $(d)$ implies $(a)$ are obvious, and $(b)$ implies $(c)$ is just Lemma \ref{nknastercapturinglemma} (which does not require $\mathfrak{m}^n_\mathcal{F}>\omega_1$). In order to prove that $(c)$ implies $(d)$, let $\mathcal{A}\in [\mathbb{P}]^{\omega_1}$. Since $\mathfrak{m}_\mathcal{F}^n>\omega_1$, so is $\mathfrak{m}(\mathbb{P})$. As $\mathbb{P}$ is $ccc$, this means that there is a filter $G$ over $\mathbb{P}$ so that $\mathcal{A}\cap G$ is uncountable. 
    \end{proof}
    \end{proposition}

\begin{corollary}Let $2\leq n\in \omega$ and $\mathcal{F}$ be an $n$-capturing construction scheme. If $\mathfrak{m}^n_\mathcal{F}>\omega_1$, then $\mathfrak{m}^n_\mathcal{F}=\mathfrak{m}_{K_m}$ for each $m\geq n$.
    
\end{corollary}
A a corollary of Theorem \ref{starfailscapturing}, we also have the following.
\begin{corollary}Let $2\leq n\in \omega$ and $\mathcal{F}$ be an $n$-capturing construction scheme. If $\mathfrak{m}^n_\mathcal{F}>\omega_1$ then $\mathfrak{m}_{K_m}=\omega_1$ for each $m<n$ and $\mathfrak{m}^m_\mathcal{F}=\omega$ for each $m>n.$
    
\end{corollary}
In order to show the consistency of $\mathfrak{m}_\mathcal{F}^n>\omega_1$, it suffices to prove that $n$-preserving $\mathcal{F}$ is a property which is preserve under finite support iteration of $ccc$ forcings.

\begin{lemma}\label{preservincapturinglemma2}Let $\mathcal{F}$ be an $n$-capturing construction scheme and $(\langle \mathbb{P}_\xi\rangle_{\xi\leq \gamma},\langle \dot{\mathbb{Q}}_\xi\rangle_{\xi<\gamma})$ be a finite support iteration of $ccc$ forcings so that $$\mathbb{P}_\xi\Vdash\text{\say{ $\dot{\mathbb{Q}}_\xi$ $n$-preserves $\mathcal{F}$ }}.$$
Then $\mathbb{P}_\gamma$ also $n$-preserves $\mathcal{F}$.
\begin{proof}The proof is carried by induction over $\gamma$ by appealing to the equivalence of $n$-capturing provided by Lemma \ref{equivalencecapturing}. Both the base and the successor steps of the induction are trivial to show. Hence, we will only do the limit case. For this, let us assume that $\gamma$ is limit and we have already showed that $\mathbb{P}_\alpha$ $n$-preserves $\mathcal{F}$ for each $\alpha<\gamma$. Let $\dot{X}$ be a $\mathbb{P}_\gamma$-name for an uncountable subset of $\omega_1$ and consider an arbitrary $q\in \mathbb{P}_\gamma$. We can take a sequence $\langle p_\xi,\alpha_\xi\rangle_{\xi\in\omega_1}\subseteq \mathbb{P}_\gamma\times \omega_1$ such that for any two distinct $\xi,\mu\in \omega_1$, the following properties hold:
\begin{itemize}
    \item $p_\xi\leq q$ and $p_\xi\Vdash\text{\say{$\alpha_\xi\in \dot{X}$ }}.$
    \item $\alpha_\xi\not=\alpha_\mu$.
\end{itemize}
By refining the sequence if necessary, we may assume that $\{ dom(p_\xi)\,:\,\xi\in\omega_1\}$ forms a $\Delta$-system with root $R\subseteq \gamma$. Now, let $\alpha<\gamma$ be such that $R\subseteq \alpha$. Observe that $q\in \mathbb{P}_\alpha$. As $\mathbb{P}_\alpha$ is $ccc$, there is $G$ a $\mathbb{P}_\alpha$-generic filter over $V$ so that $S=\{\xi\in \omega_1\,:\,p_\xi|_\alpha\in G \}$ is uncountable. In particular, this implies that for any $\xi,\mu\in S$, the conditions $p_\xi|_\alpha$ and $p_\mu|_\alpha$ are compatible. Since $R=dom(p_\xi)\cap dom(p_\mu)\subseteq \alpha$, it follows that $p_\xi$ and $p_\mu$ are compatible. Now, according to the inductive hypotheses, $\mathcal{F}$ is $n$-capturing inside $V[G]$. Therefore, there is $D\in [S]^{n-1}$ so that the family $\{\alpha_\xi\,:\,\xi\in D\}$ is captured. To finish, let $q\in \mathbb{P}_\gamma$ be such that $p\leq p_\xi$ for any $\xi\in D$. It is straightforward that $p\leq q$ and $$p\Vdash\text{\say{ $\{\alpha_\xi\,:\,\xi\in D\}\in [\dot{X}]^n$ and it is captured }}.$$

\end{proof}
\end{lemma}

\begin{theorem}\label{consistencymnfcontinuum}Let $\kappa>\omega_1$ be a regular cardinal such that $2^{<\kappa}=\kappa$. Given $\mathcal{F}$ an $n$-capturing construction scheme, there is a $ccc$-forcing $\mathbb{P}$ for which $$\mathbb{P}\Vdash\text{\say{ $\mathfrak{m}_{\mathcal{F}}^n=\mathfrak{c}=\kappa$ }}.$$
\begin{proof}We can construct a finite support iteration $(\langle \mathbb{P}_\xi\rangle_{\xi\leq \kappa},\langle \dot{\mathbb{Q}}_\xi\rangle_{\xi<\kappa})$ of $ccc$ forcings so that following properties hold for each $\xi<\kappa$:
\begin{itemize}
    \item $\mathbb{P}_\xi\Vdash \text{\say{$|\dot{\mathbb{Q}}_\xi|<\kappa$}},$
    \item $\mathbb{P}_\xi\Vdash \text{\say{ $\dot{\mathbb{Q}}_\xi$ $n$-preserves $\mathcal{F}$ }}$.
\end{itemize}
If the iteration is constructed with an appropriate bookkeeping, we can arrange that whenever $G$ is $\mathbb{P}_\kappa$-generic over $V$ and $\mathbb{Q}$ is a $ccc$ forcing inside $V[G]$ of size less than $\kappa$ which $n$-preserves $\mathcal{F}$, then there are cofinally many $\xi<\kappa$ for which $\mathbb{Q}$ is order isomorphic to $\dot{\mathbb{Q}}^G_\kappa.$ By standard arguments, all of these properties imply that $\mathbb{P}=\mathbb{P}_\kappa$ satisfies the desired conclusion.
\end{proof}
\end{theorem}

\section{Trees and gaps under $\mathfrak{m}_\mathcal{F}>\omega_1$}
\label{sectiontreesandgaps}
Recall that a partial order $(T,<)$ is called a \textit{tree} if $t_\downarrow:=\{s\in T\,:\,s<t\}$ is well-ordered for any $t\in T$. 
For each ordinal $\alpha$, the level $\alpha$ of $T$ is defined as $T_\alpha:=\{\,t\in T\,:\,ot(t_\downarrow)=\alpha\,\}.$ The \emph{height} of $T$, written as $Ht(T)$, is the minimal ordinal for which $T_\alpha=\emptyset.$ We say that $T$ is an \emph{$\omega_1$-tree} if $Ht(T)=\omega_1$,  and $T_\alpha$ is countable for each $\alpha$. An $\omega_1$-tree without uncountable chains \textit{Aronszajn} and an $\omega_1$-tree without uncountable antichains is called \emph{Suslin}.

 In \cite{treesgapsscheme}, Lopez and the third author claimed to have proved that Suslin trees exist assuming the existence of a $3$-capturing scheme. Unfortunately, the construction they provided was incorrect. In \cite{schemescruz}, we showed that Suslin trees of two special kinds, namely coherent and full, exist whenever there is a fully capturing construction scheme for which the sequence $\langle n_k\rangle_{k\in \omega}$ grows \say{quick enough} with relation to $\langle m_k\rangle_{k\in \omega}$.  We will start this section by showing that $2$-capturing schemes do not suffice to construct Suslin trees. For this task, we will need the next lemma. It is important to remark that we do not know what happens with respect to  $3$-capturing schemes.

\begin{lemma}\label{lemmacocountablecaptured}Let $\mathcal{F}$ be a $2$-capturing construction scheme. Given $X\in [\omega_1]^{\omega_1}$, the set $\{\alpha\in X\,:\, \{\beta\in X\,:\,\{\alpha,\beta\}\textit{ is captured}\}\textit{ is uncountable}\}$ is co-countable in $X.$
\begin{proof}Let $X$ be as in the hypotheses. Suppose towards a contradiction that there are uncountably many $\alpha\in X$ for which the set $C_\alpha=\{\beta\in X\,:\,\{\alpha,\beta\}\textit{ is captured }\}$ is countable. Then we can recursively construct an uncountable $Y\subseteq X$ so that if $\alpha<\beta\in Y$, then $\beta > \sup(C_\alpha\cup\{0\})$. As $\mathcal{F}$ is assumed to be $2$-capturing, there are $\alpha,\beta\in Y$ for which $\{\alpha,\beta\}$ is captured. Since $Y\subseteq X$, then $\beta\in C_\alpha$ which is a contradiction.    
\end{proof}
    
\end{lemma}

\begin{theorem}\label{nosuslinmf}Let $\mathcal{F}$ be a $2$-capturing construction scheme. If $\mathfrak{m}_\mathcal{F}>\omega_1$, then there are no Suslin trees.
\begin{proof}
Let us assume towards a contradiction that there is a Suslin tree $T.$ It is well known that Suslin trees are never Knaster. Therefore, by Propositions \ref{propequivalenceknasterpreservescheme} and \ref{lemmaequivalencefunctioncapturingpreserving}, there are $X\in [T]^{\omega_1}$ and an injective function $\zeta:X\longrightarrow \omega_1$ in such way that for any two distinct $x,y\in X$, if $\{\zeta(x),\zeta(y)\}$ is captured, then $x$ and $y$ are incompatible. We now define $$\mathbb{P}=\{ p\in[X]^{<\omega}\,:\,p\textit{ is an antichain in T}\}$$
and order it with the reverse inclusion.\\

\noindent
\underline{Claim}: $\mathbb{P}$ is $ccc$ and $2$-preserves $\mathcal{F}$.
\begin{claimproof}[Proof of claim] We will prove the claim by appealing to the Lemma \ref{lemmaequivalencefunctioncapturingpreserving}. Let $\mathcal{A}\in [\mathbb{P}]^{\omega_1}$ and $\nu:\mathcal{A}\longrightarrow \omega_1$ be an injective function. Given $p\in \mathcal{A}$, let us define $$A_p=\zeta[p]\cup \{\nu(p)\},$$ $$\alpha_p=\max(A_p).$$  Without any loss of generality we may assume that  the elements of $\{ A_p\,:\,p\in \mathcal{A}\}$ are pairwise disjoint. Furthermore, we can suppose that there are $n,k,a\in\omega$ so that the following conditions hold for any two distinct $p,q\in \mathcal{A}$:
\begin{enumerate}[label=$(\alph*)$]
    \item $|p|=n$.
    \item $\rho^{A_p}=k.$
    \item $\lVert\alpha_p\rVert_k=a.$
    \item $|A_p|=|A_q|$. Furthermore, if $h:(\alpha_p)_k\longrightarrow (\alpha_q)_k$ is the increasing bijection, then $h[A_p]=A_q$ and $h(\nu(p))=\nu(q).$
    
    \item $\rho(\alpha_p,\alpha_q)>k.$
\end{enumerate}
We now enumerate each $p\in \mathcal{A}$ as $x^p_0,\dots,x^p_{n-1}$ in such way that $\zeta(x^p_i)<\zeta(x^p_j)$ whenever $i<j.$ By virtue of the Lemma \ref{capturedfamiliestosetslemma} and the items $(b)$, $(c)$, $(d)$ and $(e)$ above, if $p,q\in \mathcal{A}$ are such that $\{\alpha_p,\alpha_q\}$ is captured, then so are $\{\nu(p),\nu(q)\}$ and $\{\zeta(x^p_i),\zeta(x^q_i)\}$. The claim follows directly from the next subclaim.\\

\noindent
\underline{Subclaim}: There are distinct $p,q\in \mathcal{A}$ so that $p\cup q$ is an antichain and $\{\alpha_p,\alpha_q\}$ are captured.
\begin{claimproof}[Proof of subclaim] Let us suppose that the subclaim is false. Then, given $p,q\in \mathcal{A}$, if $\{\alpha_p,\alpha_q\}$ is captured, then $p\cup q$ is not an antichain. In this way, there are $i_{p,q},j_{p,q}<n$ for which $x^p_{i_{p,q}}$ and $x^q_{j_{p,q}}$ are different but comparable in $T$. By means of Lemma \ref{lemmacocountablecaptured}, we can construct three sequences $\langle p_s\rangle_{s<n^2+1}\subseteq \mathcal{A}$, $\langle(i_s,j_s)\rangle_{s<n^2+1}\subseteq n\times n$ and $\langle M_s\rangle_{s<n^2+1}\subseteq [\mathcal{A}]^{\omega_1}$ so that the following properties are satisfied for each $s<n^2+1$:
\begin{itemize}
    \item For all $q\in M_s$,  $\{\alpha_{p_s},\alpha_q\}$ is captured and $(i_s,j_s)=(i_{p_s,q},j_{p_s,q}).$
    \item $p_{s+1}\in M_s$.
    \item $M_{s+1}\subseteq M_s$.
\end{itemize}
By the pigenhole principle, there are $s<r<n^2+1$ for which $(i_s,j_s)=(i_r,j_r)=(i,j)$ for some $i,j<n.$ Since $M_r$ is uncountable but ${x^{p_s}_{i}}_\downarrow\cup {x^{p_r}_{i}}_\downarrow$ is countable, there is $q\in M_r$ for which $x^{p_s}_i,x^{p_r}_i<x^q_j$. As $T$ is a tree, then $x^{p_s}_i$ and $x^{p_r}_i$ are comparable. On the other hand, $p_r\in M_s$. Thus, $\{\alpha_{p_s},\alpha_{p_r}\}$ is captured. So by the observations prior to the subclaim, this means that $\{\zeta(x^{p_s}_i),\zeta(x^{p_r}_i)\}$ is captured too. Therefore, $x^{p_s}_i$ and $x^{q_s}_i$ are incompatible. This is a contradiction, so the proof is over.
\end{claimproof}  
\end{claimproof}
 $\mathbb{P}$ is an uncountable $ccc$ forcing which $2$-preserves $\mathcal{F}$.  As $\mathfrak{m}_\mathcal{F}>\omega_1$, such forcing contains an uncountable filter, namely $G$. Note that $\bigcup G$ is an uncountable chain in $T$. This contradicts the fact that $T$ is Suslin. 
\end{proof}  
\end{theorem}

An Aronszajn tree $T$ is called \emph{special} if it can be partitioned into countably many antichains. It is a theorem of the third author (see \cite{stevotreesandlinearlyorderedsets}) that if $\mathfrak{m}_{K_2}>\omega_1$ and $T$ is an Aronszajn tree without Suslin subtrees, then $T$ is special. Since $\mathfrak{m}_\mathcal{F}\leq \mathfrak{m}_{K_2}$, we have the following corollary.
\begin{corollary}Let $\mathcal{F}$ be a $2$-capturing construction scheme. If $\mathfrak{m}_\mathcal{F}>\omega_1$, then every Aronszajn tree is special.
    
\end{corollary}
We will study $(\omega_1,\omega_1)$-gaps in the context of forcing. The reader can find fairly complete treatments of this subject in \cite{GapsandTowers}, \cite{ScheepersGaps} and \cite{StevoIlias}.

The following fundamental lemma is attributed to Kenneth Kunen. For proofs, we redirect the reader to \cite{gapsandlimits}, \cite{ScheepersGaps} and \cite{StevoIlias}.
\begin{lemma}\label{gapforcinglemma}Let $(L_\alpha,R_\alpha)_{\alpha\in\omega_1}$ be a normal $(\omega_1,\omega_1)$-pregap. Then $(L_\alpha,R_\alpha)_{\alpha\in \omega_1}$ is a gap if and only if for any $S\in [\omega_1]^{\omega_1}$ there are $\alpha<\beta\in S$ such that $(L_\alpha\cap R_\beta)\cup (L_\beta\cap R_\alpha)\not=\emptyset$.
\end{lemma}
Our next goal is to study $(\omega_1,\omega_1)$-gaps in the context of forcing.
For the rest of this section we will assume that if $(L_\alpha,R_\alpha)_{\alpha\in\omega_1}$ is an $(\omega_1,\omega_1)$-pregap then it is normal. That is, $L_\alpha\cap R_\alpha=\emptyset$ any $\alpha.$
\begin{definition}
Let $(\mathcal{L},\mathcal{R})$ be an $(\omega_1,\omega_1)$-pregap. We say that $(\mathcal{L},\mathcal{R})$ is \emph{destructible} if there is a forcing notion $\mathbb{P}$ which preserves $\omega_1$ in such way that $(\mathcal{L},\mathcal{R})$ is not a gap in some generic extension through $\mathbb{P}$. If this does not happen, the gap is said to be \emph{indestructible}.\end{definition}

Note that the Hausdorff condition is absolute. Hence, any Hausdorff gap is indestructible. In particular, this implies that under PID any $(\omega_1,\omega_1)$-gap is indestructible (see \cite{partitionpropertiesch}). 

In Definition \ref{separabilityforcingdonuts} we presented the forcing $\mathbb{P}(\mathcal{D},\mathcal{E})$ where $(\mathcal{D},\mathcal{E})$ is a pregap (not necessarily of type $(\omega_1,\omega_1)$) with $|\mathcal{D}|=|\mathcal{E}|=\omega_1$. According to Proposition \ref{separabilitygapsprop}, this forcing is $ccc$ if and only if the pregap $(\mathcal{D},\mathcal{E})$ can be separated in some $\omega_1$-preserving extension of the universe. Forcings with such properties have already been studied for the case of $(\omega_1,\omega_1)$-gaps. In the following definition we present some well-known reincarnations of them (see \cite{independenceforanalysts}, \cite{ScheepersGaps}, \cite{PartitionProblems} or \cite{yoriokadestructiblegaps}).
\begin{definition}Let $(\mathcal{L},\mathcal{R})$ be an $(\omega_1,\omega_1)$-pregap indexed as $(L_\alpha, R_\alpha)_{\alpha\in X}$ where $X$ is an uncountable set of ordinals. We define the following forcing notions:
\begin{itemize}
 \item $\chi_0(\mathcal{L},\mathcal{R})=\{p\in [X]^{<\omega}\,:\, \big(\bigcup\limits_{\alpha\in p}L_\alpha \big)\cap \big(\bigcup\limits_{\alpha\in p}R_\alpha\big)=\emptyset\}.$
    \item $\chi_1(\mathcal{L},\mathcal{R})=\{p\in [X]^{<\omega}\,:\,\forall \alpha\not=\beta\in p\,\big( (L_\alpha\cap R_\beta)\cup(L_\beta\cap R_\alpha)\not=\emptyset\big)\}.$
\end{itemize}
both ordered by reverse inclusion.
\end{definition}

The following theorem can be found in \cite{independenceforanalysts}, \cite{ScheepersGaps}, \cite{PartitionProblems} and \cite{yoriokadestructiblegaps}.

\begin{theorem}\label{cccdestructibilityequivalence}Let $(\mathcal{L},\mathcal{R})$ be an $(\omega_1,\omega_1)$-pregap:\begin{itemize}
    \item $\chi_1(\mathcal{L},\mathcal{R})$ is $ccc$ if and only if $(\mathcal{L},\mathcal{R})$ is a gap. In this case, there is some condition in $\chi_1(\mathcal{L},\mathcal{R})$ forcing $(\mathcal{L},\mathcal{R})$ to be indestructible.
    \item $\chi_0(\mathcal{L},\mathcal{R})$ is $ccc$ if and only if $(\mathcal{L},\mathcal{R})$ is destructible. In this case, there is some condition in $\chi_0(\mathcal{L},\mathcal{R})$ forcing $(\mathcal{L},\mathcal{R})$ to be separated.
\end{itemize}
\end{theorem}
A particular instance of destrucible gaps were introduced by the third author. Let $(L_\alpha,R_\alpha)_{\alpha\in\omega_1}$ be an $(\omega_1,\omega_1)$-pregap. We say that $(L_\alpha,R_\alpha)_{\alpha\in \omega_1}$ is \emph{Todor\v{c}evi\'c} if for any $S\in[\omega_1]^{\omega_1}$ there are $\alpha<\beta$ such that $L_\alpha \subseteq L_\beta$ and $R_\alpha\subseteq R_\beta$. In \cite{treesgapsscheme}, Lopez and the third author showed that the existence of Todor\v{c}evi\'c gaps follows from $CA_3$ and $CA_2(part)$. Here we give another proof of such result based on the construction given in Theorem \ref{hausdorffgapconstruction}.

\begin{theorem}[Under $CA_2(part)$]\label{todorcevicgapconstruction}Let $\mathcal{F}$ be a $\mathcal{P}$-$2$-capturing $2$-construction scheme with $\mathcal{P}=\{P_0,P_1\}$. Also, consider $(L_\alpha,R_\alpha)_{\alpha\in\omega_1}$ the Hausdorff gap constructed in Theorem \ref{hausdorffgapconstruction}. Let $$C=\bigcup\{\,\{2k,2k+1\}\,:\,k\in P_0\backslash 1\}$$ and define $$L^C_\alpha=L_\alpha\cap C=\{ 2k+\Xi_\alpha(k)\,:\,k\in P_0\backslash 1\textit{ and }\Xi_\alpha(k)\geq 0\},$$
$$R^C_\alpha=R_\alpha\cap C=\{2k+(1-\Xi_\alpha(k))\,:\,k\in P_0\backslash 1,\textit{ and }\Xi_\alpha(k)\geq 0\}$$
for each $\alpha\in \omega_1$. Then $(\mathcal{L}|_C,\mathcal{R}|_C)=(L^C_\alpha,R^C_\alpha)_{\alpha\in \omega_1}$ is a Todor\v{c}evi\'c gap.
\begin{proof}$(L^C_\alpha,R^C_\alpha)_{\alpha\in\omega_1}$ is a pregap by similar reasons as the ones in the proof of Theorem \ref{hausdorffgapconstruction}. That is, $L^C_\alpha\cap R^C_\alpha=\emptyset$ for each $\alpha\in\omega_1$ and if $\alpha<\beta$ then:
\begin{enumerate}[label=$(\alph*)$] 

\item  $L^C_\alpha\backslash L^C_\beta\subseteq \{2k+\Xi_\alpha(k)\,:\,k\leq \rho(\alpha,\beta)\textit{ and }\Xi_\alpha(k)\not=\Xi_\beta(k)\},$
\item $R^C_\alpha\backslash R^C_\beta\subseteq \{ 2k+(1-\Xi_\alpha(k))\,:\,k\leq \rho(\alpha,\beta)\textit{ and }\Xi_\alpha(k)\not=\Xi_\beta(k)\}.$ 
\end{enumerate}
In particular, as a consequence of (a) we have:
\begin{enumerate}
\item[$(c)$] $L^C_\alpha\cap R^C_\beta\subseteq\{2k+\Xi_\alpha(k)\,:\,k\leq \rho(\alpha,\beta)\textit{ and }\Xi_\alpha(k)\not=\Xi_\beta(k)\}.$

\end{enumerate}

In order to prove that $(L^C_\alpha,R^C_\alpha)_{\alpha\in \omega_1}$ is a gap we will use the equivalence provided by Lemma \ref{gapforcinglemma}. Let $S\in [\omega_1]^{\omega_1}$. Since $\mathcal{F}$ is $\mathcal{P}$-2-capturing there are $\alpha<\beta\in S$ so that $\{\alpha,\beta\}$ is captured at some level $l\in P_0$.  It follows that $\Delta(\alpha,\beta)=l=\rho(\alpha,\beta)$. Hence, $\Xi_\alpha(k)=\Xi_\beta(k)$ for any $k<l$. From this and by the the point $(c)$ above we deduce that $L^C_\alpha \cap R^C_\beta\subseteq \{2l+\Xi_\alpha(l)\}$. Furthermore, since $l=\rho(\alpha,\beta)$ and $\mathcal{F}$ is a construction scheme, we also have that $\Xi_\alpha(l)=0$ and $\Xi_\beta(l)=1$. Therefore we conclude that $L^C_\alpha\cap R^C_\beta=\{2l\}$, so in particular $(L^C_\alpha\cap R^C_\beta)\cup (L^C_\beta\cap R^C_\alpha)\not=\emptyset.$ As we said before, by Lemma \ref{gapforcinglemma} we are done.

Now we will prove that $(L^C_\alpha,R^C_\alpha)_{\alpha\in\omega_1}$ is Todor\v{c}evi\'c. Let $S\in [\omega_1]^{\omega_1}$. Since $\mathcal{F}$ is $\mathcal{P}$-2-capturing, there are $\alpha<\beta\in S$ such that $\{\alpha,\beta\}$ is captured at some level $l\in P_1$. Again, by the points (a) and (b) written at the beginning of the proof we deduce that $L^C_\alpha\backslash L^C_\beta\subseteq \{2l\} $ and $R^C_\alpha\backslash R^C_\beta\subseteq \{2l+1\}$. But $l\in P_1$ so $\{2l,2l+1\}$ has empty intersection with both $L^C_\alpha$ and $R^C_\alpha$ by definition of these two sets. In this way we conclude that $L^C_\alpha\backslash L^C_\beta=\emptyset$ and $R^C_\alpha\backslash R^C_\beta=\emptyset.$ In other words, $L^C_\alpha\subseteq L^C_\beta$ and $R^C_\alpha\subseteq R^C_\beta$. 
\end{proof}
\end{theorem}
It is natural to ask whether destructible gaps can be constructed only by using a $2$-capturing construction scheme. In the next theorem we will show that this is not the case. 
In order to prove it, we will make use of the following two lemmas which are implicit in the proof of Theorem \ref{cccdestructibilityequivalence}.

\begin{lemma}\label{lemmachi0mf}Let $(\mathcal{L},\mathcal{R})=(L_\alpha,R_\alpha)_{\alpha\in\omega_1}$ be an $(\omega_1,\omega_1)$-gap. Suppose that $\mathcal{A}\in[\chi_0(\mathcal{L},\mathcal{R})]^{\omega_1}$. Then there are $n\in\omega$ and  $\mathcal{B}\in [\mathcal{A}]^{\omega_1}\cap \mathscr{P}([\omega_1]^n)$ an uncountable root-tail-tail $\Delta$-system such that for any two distinct $p,q\in  \mathcal{B}$, the following conditions are equivalent:
\begin{itemize}
    \item $p$ is compatible with $q$.
    \item $(L_{p(n-1)}\cap R_{q(n-1)})\cup (L_{q(n-1)}\cap R_{p(n-1)})=\emptyset.$ In other words, $\{p(n-1)\}$ is compatible with $\{q(n-1)\}.$
\end{itemize}
\end{lemma}

\begin{lemma}\label{lemmachi1mf}Let $(\mathcal{L},\mathcal{R})=(L_\alpha,R_\alpha)_{\alpha\in\omega_1}$ be an $(\omega_1,\omega_1)$-gap. Suppose that $\mathcal{A}\in[\chi_1(\mathcal{L},\mathcal{R})]^{\omega_1}$. Then there are $n,r\in\omega$ and  $\mathcal{B}\in [\mathcal{A}]^{\omega_1}\cap \mathscr{P}([\omega_1]^n)$  a root-tail-tail $\Delta$-system with root $R$ of cardinality $r$ such that for any two distinct $p,q\in  \mathcal{B}$, the following conditions are equivalent:
\begin{itemize}
    \item $p$ is compatible with $q$.
    \item $(L_{p(r)}\cap R_{q(r)})\cup (L_{q(r)}\cap R_{p(r)})\not=\emptyset.$ In other words, $\{p(r)\}$ is compatible with $\{q(r)\}.$ 
\end{itemize}
    
\end{lemma}
\begin{theorem}\label{nodestructiblemf}Let $\mathcal{F}$ be a $2$-capturing construction scheme. If $\mathfrak{m}_\mathcal{F}>\omega_1$, then there are no destructible gaps.
\begin{proof}Let us assume towards a contradiction that there is a destructible gap $(\mathcal{L},\mathcal{R})=(L_\alpha,R_\alpha)_{\alpha\in\omega_1}.$ The forcing $\chi_0(\mathcal{L},\mathcal{R})$ is never Knaster. Therefore, by Propositions \ref{propequivalenceknasterpreservescheme} and \ref{lemmaequivalencefunctioncapturingpreserving}, there are $X'\in [\chi_0(\mathcal{L},\mathcal{R})]^{\omega_1}$ and an injective function $\zeta':X\longrightarrow \omega_1$ in such way that for any two distinct $p,q\in X'$, if $\{\zeta'(p),\zeta'(q)\}$ is captured, then $p$ and $q$ are incompatible in $\chi_0(\mathcal{L},\mathcal{R})$. By virtue of the Lemma \ref{lemmachi0mf}, we may assume that the elements $X'$ are singletons. In this way, we can define $X=\{\alpha\in \omega_1\,:\,\{\alpha\}\in X'\}$ and $\zeta:X\longrightarrow \omega_1$ given by $\zeta(\alpha)=\zeta'(\{\alpha\})$.
Note that for any two distinct $\alpha,\beta\in X$, if $\{\zeta(\alpha),\zeta(\beta)\}$ is captured, then $$(L_\alpha\cap R_\beta)\cup (L_\alpha\cap R_\beta)\not=\emptyset.$$
\noindent
\underline{Claim}: $\chi_1( L_\alpha,R_\alpha)_{\alpha\in X}$ is $ccc$ and $n$-preserves $\mathcal{F}.$
\begin{claimproof}[Proof of claim]We will prove the claim by appealing to the Lemma  \ref{lemmaequivalencefunctioncapturingpreserving}. Let $\mathcal{A}$ be an uncountable subset of $\chi_1(L_\alpha,R_\alpha)_{\alpha\in X}$ and  $\nu:\mathcal{A}\longrightarrow \omega_1$ be an injective function. According to the Lemma \ref{lemmachi1mf}, we may suppose that there are $n,r\in \omega$ so that $\mathcal{A}\in [X]^n$ and it forms a root-tail-tail $\Delta$-system with root $R$ of cardinality $r$ so that for any two distinct $p,q\in \mathcal{A}$, if $(L_{p(r)}\cap R_{q(r)})\cup (L_{q(r)}\cap R_{p(r)})\not=\emptyset$, then $p$ and $q$ are compatible.  Given $p\in \mathcal{A}$, let us define $$A_p=\ \{\zeta(p(r)),\nu(p)\},$$ $$\alpha_p=\max(A_p).$$ 
Since both $\zeta$ and $\nu$ are injective functions and $p(r)\not=q(r)$ whenever $p\not=q$, we may assume that the $A_p$'s are pairwise disjoint. Even more,we can suppose that there are $k,a\in\omega$ so that the following conditions hold for any two distinct $p,q\in \mathcal{A}$:
\begin{enumerate}[label=$(\alph*)$]
    \item $\rho^{A_p}=k.$
    \item $\lVert\alpha_p\rVert_k=a.$
    \item If $h:(\alpha_p)_k\longrightarrow (\alpha_q)_k$ is the increasing bijection $h(\nu(p))=\nu(q)$ and $h(\zeta(p(r)))=\zeta(q(r))$.
      \item $\rho(\alpha_p,\alpha_q)>k.$
\end{enumerate}
Since $\mathcal{F}$ is assumed to be $2$-capturing, there are distinct $p_0,p_1\in \mathcal{A}$ for which $\{\alpha_{p_0},\alpha_{p_1}\}$ is captured. By virtue of the Lemma \ref{capturedfamiliestosetslemma} and the points $(a)$, $(b)$, $(c)$ and $(d)$ above, are $\{\nu(p_0),\nu(p_1)\}$ and $\{\zeta(p_0(r)),\zeta(p_1(r))\}$ are also captured. In particular, this means that $(L_{p_0(r)}\cap R_{p_1(r)})\cup (L_{p_1(r)}\cap R_{p_0(r)})\not=\emptyset.$ Thus, $p_1$ and $p_1$ are compatible. This finishes the proof. 

\end{claimproof}
By the previous claim $\chi_1(L_\alpha,R_\alpha)_{\alpha\in X}$ is an uncountable $ccc$ forcing which $2$-preserves $\mathcal{F}$. Since $\mathfrak{m}_\mathcal{F}>\omega_1$, it follows that $(L_\alpha,R_\alpha)_{\alpha\in X}$ is an indestructible gap due to the Lemma \ref{cccdestructibilityequivalence} (hence, so is $(\mathcal{L},\mathcal{R})$). This is a contradiction, so the proof is over.
\end{proof}
\end{theorem}

\section{Capturing with partitions}\label{capturingwithpartitionssection}

Our next goal is to show construction schemes which are $n$-capturing may not be $\mathcal{P}$-$n$-capturing for any non-trivial partition $\mathcal{P}$ of $\omega$. For this, we introduce a natural filter over $\omega$ which is definable from a construction scheme.
\begin{definition}Let $2\leq n \in \omega$ and $\mathcal{F}$ be a construction scheme. Given $S\subseteq \omega_1$, we define the \textit{$n$-projection} of $S$ as $$\pi_n(S)=\{ \rho^D\,:\,D\in [S]^n\textit{ and }D\textit{ is captured }\}$$
In other words, $\pi_n(S)$     is the set of all $l\in\omega$ for which there is $D\in [S]^n$ which is captured at level $l.$
\end{definition}
\begin{rem}Given $S\subseteq \omega_1 $, the $n$-projection of $S$ can be recovered from any cofinal subset of $[S]^{<\omega}$. That is, if $G$ is a cofinal subset of $[S]^{<\omega}$, then  $$\pi_n(S)=\bigcup_{D\in G}\pi_n(D).$$
    
\end{rem}

\begin{definition}[The $n$-projection filter]Let  $n\in \omega$ and $\mathcal{F}$ be a construction scheme. We define $\mathcal{U}_n(\mathcal{F})$ as the set of all $A\subseteq \omega$ for which there is $S\in [\omega_1]^{\omega_1}$ such that $$\pi_n(S)\subseteq A.$$ 
\end{definition}
\begin{rem}Note that if $\mathcal{F}$ is not $n$-capturing then $\emptyset\in \mathcal{U}_n(\mathcal{F})$. Therefore, $\mathcal{U}_n(\mathcal{F})$ is not actually a filter in this case. 
\end{rem}

\begin{proposition}\label{lemmafilteruf}Let $2\leq n \in \omega$ and $\mathcal{F}$ be a construction scheme. If $\mathcal{F}$ is $n$-capturing, then $\mathcal{U}_n(\mathcal{F})$ is a non-principal filter over $\omega$.
\begin{proof}Since $\mathcal{F}$ is $n$-capturing, it follows that each member of $\mathcal{U}_n(\mathcal{F})$ is non-empty. Thus, in order to prove that $\mathcal{U}_n(\mathcal{F})$ is a filter, it is enough to show that the family $\{\,\pi_n(S)\,:\,S\in[\omega_1]^{\omega_1}\,\}$ is downwards directed with respect to $\subseteq$. Let $S,S'\in [\omega_1]^{\omega_1}$. We will prove that there is $A\in [\omega_1]^{\omega_1}$ such that $\pi_n(A)\subseteq \pi_n(S)\cap \pi_n(S')$. For this, first note that we can recursively construct a sequence $\langle \alpha_\xi,\beta_\xi\rangle\subseteq S\times S'$  so that  $\alpha_\xi<\beta_\xi<\alpha_\zeta$
for each $\xi<\zeta\in \omega_1$. By refining the sequence if necessary, we may assume without loss of generality that there is $k\in\omega$ so that for any two distinct $\xi,\mu\in \omega_1$, the following conditions hold:
\begin{enumerate}[label=$(\arabic*)$]
    \item $\rho(\alpha_\xi,\beta_\xi)<k,$
    \item $\lVert \alpha_\xi\rVert_k=\lVert \alpha_\mu\rVert_k$ and $\lVert \beta_\xi\rVert_k=\lVert \beta_\mu\rVert_k$. In particular, $\rho(\beta_\mu,\beta_\xi)>k$.
\end{enumerate}
Let $A=\{\beta_\xi\,:\,\xi\in \omega_1\}$. Then $A\in[ S']^{\omega_1}$, so $\pi_n(A)\subseteq \pi_n(S')$. We claim that $\pi_n(A)\subseteq \pi_n(S)$. For this purpose, let $l\in \pi_n(A)$ and consider $D\in[A]^n$ so that $\{ \beta_\xi\,:\xi\in D\}$ is captured at level $l$. Note that  $l>k$ due to the condition (2).\\
 
\noindent
\underline{Claim}: If $\xi,\mu\in D$ are distinct then $\Delta(\alpha_\xi,\alpha_\mu)=l=\rho(\alpha_\xi,\alpha_\mu)$.
\begin{claimproof}[Proof of claim] Let $h:(\beta_\xi)_k\longrightarrow (\beta_\mu)_k$ be the increasing bijection. By the point (1),  we have that $\alpha_\xi\in (\beta_\xi)_k$ and $\alpha_\mu\in (\beta_\mu)_k$. Furthermore, $h(\alpha_\xi)=\alpha_\mu$ by the point (2).  Since $\alpha_\xi\not=\alpha_\mu$, we may use Lemma \ref{lemmahdeltarhoinequalities} to conclude that $$l=\rho(\beta_\xi,\beta_\mu)\geq \rho(\alpha_\xi,\alpha_\mu)\geq \Delta(\alpha_\xi,\alpha_\mu)\geq \Delta(\beta_\xi,\beta_\mu)=l.$$
So we are done.
\end{claimproof}
Note that $\Xi_{\alpha_\xi}(l)\geq 0$ for each $\xi\in D$. This is due to the point (b) of Lemma \ref{lemmaxi}. By the point (c)  of such lemma, $\Xi_{\alpha_\xi}(l)=\Xi_{\beta_\xi}(l)$. Thus, $\{\alpha_\xi\,:\xi\in D\}$ is captured at level $l$. That is, $l\in \pi_n(S)$. As $l\in \pi_n(A)$ was arbitrary, we get that  $\pi_n(A)\subseteq \pi_n(S)$.\\

 Now, we will show that $\mathcal{U}_n(\mathcal{F})$ is non-principal. Let $S\in [\omega_1]^{\omega_1}$ and  $k\in \omega$. Then there is $S'\in [S]^{\omega_1}$ such that $\lVert \alpha\rVert_k=\lVert \beta\rVert_k$ for all $\alpha,\beta\in S$. This implies that $\rho^D>k$ for any $D\in [S']^n$. Therefore, $\pi_n[S']\subseteq \pi_n[S]\backslash k$. This finishes the proof.
\end{proof}   
\end{proposition}
\begin{rem}\label{remarkpartitionpositive}Suppose that $\mathcal{F}$ is an $n$-capturing construction scheme. Note that $A\in \mathcal{U}_n(\mathcal{F})^+$ if and only if for any $S\in [\omega_1]^{\omega_1}$ there are infinitely many $l\in A$ for which there is $D\in [S]^n$ which is captured at level $l$. In this way,  $\mathcal{F}$ is $\mathcal{P}$-$n$-capturing for a partition $\mathcal{P}$ of $\omega$ if and only if $\mathcal{P}\subseteq \mathcal{U}_n(\mathcal{F})^+$.
\end{rem}
As a direct corollary of the previous remark we have that:
\begin{corollary}Let $2\leq n\in \omega$ and $\mathcal{F}$ be a $n$-capturing construction scheme. There is a non-trivial partition $\mathcal{P}$ of $\omega$ for which $\mathcal{F}$ is $\mathcal{P}$-$n$-capturing if and only if $\mathcal{U}_n(\mathcal{F})$ is not an ultrafilter.
\end{corollary}
Now, we will show that it is consistent the existence of a construction scheme which is $n$-capturing and $\mathcal{P}$-$(n-1)$-capturing for some non-trivial partition of $\omega$, but it is not $\mathcal{P}'$-$n$-capturing for any non-trivial partition $\mathcal{P}'$.  This will be done by  starting with a $\mathcal{P}$-$n$-capturing construction scheme, and building a suitable finite support iteration of $ccc$ forcings which force $\mathcal{U}_{n}(\mathcal{F})$ to be an ultrafilter. In the following definition, we describe the forcings that will be used for this task.

\begin{definition}Let $\mathcal{F}$ be a construction scheme. Given $2\leq n\in \omega$ and $A\subseteq \omega$, we define the forcing $\mathbb{D}_n(\mathcal{F},A)$ as the family of all $p\in \text[\omega_1]^{<\omega}$ with the following property:
\begin{center}
There is no $D\in[p]^n$ such that $D$ is captured at some level $l\in A$.
\end{center}
We order $\mathbb{D}_n(\mathcal{F},A)$ with respect to $\subseteq$.
\end{definition}
The previous forcing was considered in \cite{forcingandconstructionschemes} for the particular case where $A=\omega$.

\begin{lemma}\label{Dnccc}Let  $2\leq n\in \omega$ and $\mathcal{F}$ be a $2$-capturing construction scheme. Then $\mathbb{D}_n(\mathcal{F},A)$ is $ccc$ for any $A\subseteq \omega$.
\begin{proof}Let $\mathcal{A}$ be an uncountable subset of $\mathbb{D}_n(\mathcal{F},A)$. Consider  $\mathcal{A}'\subseteq \mathcal{A}$ an uncountable root-tail-tail $\Delta$-system with root $R$ so that any two elements of $\mathcal{A}'$ have the same cardinality. Note that  we can enumerate $\mathcal{A}'$ as $\langle p_\alpha\rangle_{\alpha\in \omega_1}$ in such way that $\max(p_\alpha\backslash R)<\min(p_\beta\backslash R)$ whenever $\alpha<\beta.$ Since $\mathcal{F}$ is $2$-capturing, we can find $\delta<\gamma\in \text{Lim}$ so that the set $\{p_\delta\cup p_{\delta+1},p_\gamma\cup p_{\gamma+1}\}$ is captured at some level $l$. The proof will end by showing the following claim.\\

\noindent
\underline{Claim:} $p=p_\delta\cup p_{\gamma+1}\in \mathbb{D}_n(\mathcal{F},A)$.
\begin{claimproof}[Proof of claim] Let $D\in [p]^n$ which is captured and $h:(p_\delta\cup p_{\delta+1})_{l-1}\longrightarrow (p_\gamma\cup p_{\gamma+1})_{l-1}$ be the increasing bijection. It is easy to see that $h[p_\delta]=p_\gamma$ and $h[p_{\delta+1}]=p_{\gamma+1}$. Therefore, $$\lVert \alpha\rVert_{l-1}=\lVert h(\alpha)\rVert_{l-1}<\lVert \beta\rVert_{l-1} $$
for each $\alpha\in p_\delta\backslash R$ and $\beta\in p_{\gamma+1}\backslash R.$ Thus, $\rho^D<l$. In this way,  either $D\subseteq p_\delta$ or $D\subseteq p_{\gamma+1}$. As both $p_\delta$ and $p_{\gamma+1}$ belong to $\mathbb{D}_n(\mathcal{F},A)$, it follows that $\rho^D\notin A$. 
 \end{claimproof}
\end{proof} 
\end{lemma}
The following corollary will help us turn $\mathcal{U}_n(\mathcal{F})$ into an ultrafilter.
\begin{corollary}\label{notpositiveforcinglemma} Let $2\leq n\in \omega$, $\mathcal{F}$ be a $2$-capturing construction scheme and $A\subseteq \omega$. Then there is a condition $p\in \mathbb{D}_n(\mathcal{F},A)$ so that $$p\Vdash \text{\say{ $A\not\in \mathcal{U}_n(\mathcal{F})^+$ }}.$$
\begin{proof}Since $\mathcal{F}$ is $2$-capturing, $\mathbb{D}_n(\mathcal{F},A)$ is a $ccc$-forcing. Furthermore, it is uncountable. Thus, there is $p\in \mathbb{D}_n(\mathcal{F},A)$ which forces the generic filter to be uncountable. It follows that if $G$ is a generic filter over $V$ so that $p\in G$, then $S_G=\bigcup G$ is an uncountable subset of $\omega_1$ with $\pi_n(S_G)\cap A=\emptyset.$ Thus, in $V[G]$, $A\notin \mathcal{U}_{n+1}(\mathcal{F})$.
\end{proof}
\end{corollary}

\begin{lemma}\label{Dncapturingpreservinglemma}Let $2\leq n\in \omega$ and $\mathcal{F}$ be an $n$-capturing construction scheme. If $A\in \mathcal{U}_n(\mathcal{F})^+$ then $\mathbb{D}_n(\mathcal{F},\omega\backslash A)$ $n$-preserves $\mathcal{F}.$
\begin{proof}We will prove this lemma by appealing to the equivalence provided by Lemma \ref{lemmaequivalencefunctioncapturingpreserving}. Let $\mathcal{A}$ be an uncountable subset of $\mathbb{D}_n(\mathcal{F},\omega\backslash A)$ and $\nu:\mathcal{A}\longrightarrow \omega_1$ be an injective function.
Given $p\in \mathcal{A}$, let $D_p=p\cup \{\nu(p)\}$ and $\alpha_p=\max(D_p)$. By refining $\mathcal{A}$ if necessary, we may assume that are $j\in \omega$, $a<m_j$, $C\subseteq a+1$ and $b\in C$ so that the following conditions hold for any $p\in \mathcal{A}$:
\begin{enumerate}[label=$(\alph*)$]
    \item $\rho^{D_p}\leq j,$
    \item $\lVert \alpha_p\rVert_j=a$,
    \item $(\alpha_p)_j[C]=p$ and $(\alpha_p)_j(b)=\nu(p)$.
\end{enumerate}
As $A\in \mathcal{U}_n(\mathcal{F})^+$, there are $l\in A$ and $\{p_0,\dots,p_{n-1}\}\in [\mathcal{A}]^n$ for which $\{\,\alpha_{p_i}\,:i<n\,\}$ is captured at level $l$. Note that $l>j$ due to the point (b) above. According to Lemma \ref{capturedfamiliestosetslemma}, the family $\{D_{p_i}\,:\,i<n\}$ is  captured at level $l$. From this fact and by the point (c), the same holds for both $\{p_0,\dots,p_{n-1}\}$ and $\{\nu(p_0),\dots,\nu(p_{n-1})\}$. The proof follows from the following claim. \\

\noindent
\underline{Claim}: $p=\bigcup\limits_{i<n}p_i\in \mathbb{D}_n(\mathcal{F},\omega\backslash A).$
\begin{claimproof}[Proof of claim]Let  $D\in [p]^n$ which is captured and consider $F\in \mathcal{F}_l$ such that $p\subseteq F$. As $D\subseteq p$, then $\rho^D\leq \rho^p=l$. If $\rho^D=l$, we are done because $l\not\in \omega\backslash A.$ On the other hand, if  $\rho^D<l$ then $D\subseteq F_{\Xi_D(l)}$.  Since  the family $\{p_0,\dots,p_n\}$ is captured at level $l$, there is $i<n$ so that   $F_{\Xi_D(l)}\cap p=p_i$. In this way, $D\subseteq p_i$. Thus, as $p_i\in \mathbb{D}_n(\mathcal{F},\omega\backslash A)$ then $\rho^D\not\in \omega\backslash A$.
\end{claimproof}
\end{proof}
\end{lemma}

\begin{lemma}\label{preservingpositivelemma1}Let $3\leq n\in \omega$, $\mathcal{F}$ be an $n$-capturing construction scheme and $A \subseteq \omega$. If $B\in \mathcal{U}_{n-1}(\mathcal{F})^+$, then  $$\mathbb{D}_n(\mathcal{F},A)\Vdash \text{\say{ $B\in \mathcal{U}_{n-1}(\mathcal{F})^+$ }}.$$

\begin{proof}Let $\dot{X}$ be a name for an  uncountable subset of $\omega_1$ and $q\in \mathbb{D}_{n}(\mathcal{F}, A)$. Using that $B\in \mathcal{U}_{n-1}(\mathcal{F})^+$ and by arguing in a similar way as in the previous lemmas, we can find a sequence $\langle p_i,\alpha_i\rangle_{i<{n-1}}\subseteq \mathbb{D}_n(\mathcal{F},A)\times \omega_1$ so that:
\begin{enumerate}[label=$(\arabic*)$]
    \item both $\{p_0,\dots,p_{n-2}\}$ and $\{\alpha_0,\dots,\alpha_{n-2}\}$ are captured at some level $l\in B$,
    \item $\alpha_i\not=\alpha_j$ whenever $i\not=j$,
    \item $p_i\leq q$ and $p_i\Vdash\text{\say{ $\alpha_i\in \dot{X}$ }}$ for each $i<n-1.$ 
    \end{enumerate}
    It is easy to see that $p=\bigcup\limits_{i<n-1}p_i$ is a condition in the forcing that we are considering. Furthermore, $p\leq q$ and $p\Vdash \text{\say{ $\{\alpha_0,\dots,\alpha_{n-1}\}\in [\dot{X}]^{n-1}$ }}.$ This means that $p\Vdash \text{\say{ $\pi_{n-1}(\dot{X})\cap B\not=\emptyset$ }}.$ Thus, the proof is over.
\end{proof}
\end{lemma}

The following lemma is proved in the exact same way as the previous one. For that reason, we leave the proof to the reader.
\begin{lemma}\label{preservingpositiveslemma2}Let $\mathcal{F}$ be an $n$-capturing construction scheme, $B\in \mathcal{U}_n(\mathcal{F})^+$ and $(\langle \mathbb{P}_\xi\rangle_{\xi\leq \gamma},\langle \dot{\mathbb{Q}}_\xi\rangle_{\xi<\gamma})$ be a finite support iteration of $ccc$ forcings so that $$\mathbb{P}_\xi\Vdash\text{\say{ $\dot{\mathbb{Q}}_\xi\Vdash\text{\say{$B\in \mathcal{U}_{n}(\mathcal{F})^+$}}$}}.$$  Then $\mathbb{P}_\gamma\Vdash\text{\say{$B\in \mathcal{U}_{n}(\mathcal{F})^+$}}$.
\end{lemma}

By combining all the results we have so far, we get the following theorem.
\begin{theorem}Let $\mathcal{F}$ be $\mathcal{P}$-$n$-capturing construction scheme. There is a $ccc$ forcing $\mathbb{P}$ satisfies the following properties: $$\mathbb{P}\Vdash\text{\say{$\mathcal{F}$ is $n$-capturing and $\mathcal{P}$-$(n-1)$-capturing}},$$
$$\mathbb{P}\Vdash\text{\say{$\mathcal{U}_n(\mathcal{F})$ is an ultrafilter }}.$$
In particular, $\mathbb{P}$ forces that $\mathcal{F}$ is not $\mathcal{P}'$-$n$-capturing for any non-trivial partition $\mathcal{P}'$ of $\omega$.
\begin{proof}Let $\kappa=\mathfrak{c}$. We can construct a finite support iteration  $(\langle \mathbb{P}_\xi\rangle_{\xi\leq \kappa},\langle \dot{\mathbb{Q}}_\xi\rangle_{\xi<\kappa})$  of forcings so that given $\xi<\kappa$, the following condition holds:
$$\mathbb{P}_\xi\Vdash\text{\say{$\dot{\mathbb{Q}_\xi}=\mathbb{D}_n(\mathcal{F},\dot{A})$ for some $\dot{A}\in[\omega]^{\omega}$ with $\dot{A},\omega\backslash \dot{A}\in \mathcal{U}_n(\mathcal{F})^+$}}.$$ 
By means of Lemma \ref{Dnccc}, it follows that $\mathbb{P}_\kappa$ is $ccc$. According to Lemma \ref{Dncapturingpreservinglemma} and Lemma \ref{preservincapturinglemma2}, we have that $\mathbb{P}$ $n$-preserves $\mathcal{F}$. Furthermore, $\mathbb{P}\Vdash\text{\say{ $\mathcal{P}\subseteq\mathcal{U}_{n-1}(\mathcal{F})^+$}}$ due to Lemma \ref{preservingpositivelemma1} and Lemma\ref{preservingpositiveslemma2}. If we construct the iteration with an appropriate bookkeeping and make use of Corollary \ref{notpositiveforcinglemma}, we can arrange $\mathbb{P}_\kappa$ to force   $\mathcal{U}_n(\mathcal{F})$ to be an ultrafilter.
This finishes the proof.
\end{proof}
\end{theorem}
\begin{corollary}\label{partitionncapturingcorolary}Let $1\leq n\in \omega$. It is consistent that there is a construction scheme $\mathcal{F}$ such that:\begin{itemize}
    \item There is a non-trivial partition of $\omega$, namely $\mathcal{P}$, such that $\mathcal{F}$ is $\mathcal{P}$-$n$-capturing.
    \item $\mathcal{F}$ is $n+1$-capturing but it is not $\mathcal{P}'$-$n+1$-capturing for any non-trivial partition $\mathcal{P}'$.
\end{itemize}
\end{corollary}

    \section{Ramsey ultrafilters from construction schemes}\label{sectionramseyultafilters}
Let $\mathcal{U}$ be an ultrafilter over $\omega$. Recall that $\mathcal{U}$ is \textit{Ramsey} if for each infinite partition  $\mathcal{P}\subseteq \mathcal{U}^*$ of $\omega$, there is $A\in \mathcal{U}$ so that $|A\cap P|\leq 1$
    for all $P\in \mathcal{P}.$ Here, $\mathcal{U}^*$ denotes the dual ideal of $\mathcal{U}$. That is, the set $\{A\in \omega\,:\,\omega\backslash A\in \mathcal{}U\}$. Ramsey ultrafilters, also called selective, play a highly important roll in modern set theory. In \cite{mapsontrees}, the third author presented clever ways of defining filters over $\omega_1$ by using walks on ordinals. In \cite{walksultrafilters}, he used such filters in order to show that the existence of Ramsey ultrafilters follows from $\mathfrak{m}>\omega_1$. This result is quite interesting since the existence of ultrafilters with particular properties usually follows from equalities of the form \say{$\textit{cardinal invariant}=\mathfrak{c}$} instead of inequalities of the form \say{$\textit{cardinal invariant}>\omega_1$}. In this section we will show that given a $2$-capturing construction scheme $\mathcal{F}$, the filter $\mathcal{U}_n(\mathcal{F})$ is a Ramsey ultrafilter under $\mathfrak{m}^n_\mathcal{F}>\omega_1$. Furthermore, we will define the filter $\mathcal{V}(\mathcal{F})$ and show that, under $\mathfrak{m}_\mathcal{F}>\omega_1$, this filter is an ultrafilter over $\omega_1$ which projects into a Ramsey ultrafilter. Beyond the hypothesis imposed over the cardinal invariant $\mathfrak{m}^n_\mathcal{F}$, the most important feature of the constructions below is that both ultrafilters are explicitly definable using a combinatorial structure over $\omega_1$.

    \begin{theorem}\label{UnFramseytheorem}Let $2\leq n\in \omega$ and $\mathcal{F}$ be an $n$-capturing construction scheme for which $\mathfrak{m}^n_\mathcal{F}>\omega_1$. Then $\mathcal{U}_n(\mathcal{F})$ is a Ramsey ultrafilter.
    \begin{proof} Suppose towards a contradiction that $\mathcal{U}_n(\mathcal{F})$ is not an ultrafilter and let $A\subseteq \omega$ be such that $A, \omega\backslash A\in \mathcal{U}_n(\mathcal{F})^+$. According to the Lemmas \ref{Dnccc} and \ref{Dncapturingpreservinglemma}, $\mathbb{D}_n(\mathcal{F},A)$ is a $ccc$ forcing so that $\mathbb{D}_n(\mathcal{F},A)\Vdash\text{\say{ $\mathcal{F}$ is $n$-capturing }}$. Hence, its Martin's number is greater than $\omega_1$ due to the hypotheses. In this way, there is an uncountable filter $G$ over $\mathbb{D}_n(\mathcal{F},A)$. Let $S_G=\bigcup G$. Then $S_G$ is an uncountable subset of $\omega_1$ satisfying that $\pi_n(S_G)\cap A=\emptyset.$ This contradicts the fact that $A\in \mathcal{U}_n(\mathcal{F})^+$.

    Now, we will show that $\mathcal{U}_n(\mathcal{F})$ is Ramsey. For this purpose, let $\langle P_i\rangle_{i\in \omega}$ be a partition of $\omega$ so that $P_i\notin \mathcal{U}_n(\mathcal{F})$ for every $i\in \omega$. It suffices to prove that there is $S\in [\omega_1]^{\omega_1}$ with $|\pi_n(S)\cap P_i|\leq 1$ for any $i\in \omega$. Let $\mathbb{P}$ be the forcing consisting of all $p\in[\omega_1]^{<\omega}$ so that: $$|\pi_n(p)\cap A_i|\leq 1$$
    for each $i\in \omega$. It is easy to see that if $G$ is an uncountable filter over $\mathbb{P}$ then  $S_G=\bigcup G$ satisfies the required property.  As we are assuming that $\mathfrak{m}_\mathcal{F}^n>\omega_1$, the existence of such $G$ will follow if $\mathbb{P}$ is $ccc$ and $\mathbb{P}\Vdash\text{\say{ $\mathcal{F}$ is $n$-capturing }}$. The proof this two facts is similar to ones in Lemmas \ref{Dnccc} and \ref{Dncapturingpreservinglemma}. For that reason, we leave the details to the reader.   
    \end{proof} 
    \end{theorem}
    \begin{definition}[A square bracket operation]Let $\mathcal{F}$ be a construction scheme. We define $\llbracket\cdot,\cdot\rrbracket:[\omega_1]^2\longrightarrow \omega$ as:
    $$\llbracket \alpha,\beta\rrbracket=\min(\,(\beta)_{\rho(\alpha,\beta)-1}\backslash \alpha\,)$$
    for each $\alpha<\beta$. Given $S\subseteq \omega_1$, we also define $$\llparenthesis S\rrparenthesis=\{\,\llbracket \alpha,\beta\rrbracket \,:\,\alpha,\beta\in S\textit{ and }\{\alpha,\beta\} \textit{ is captured }\}$$
    \end{definition}
\begin{rem}\label{remarksqbracket1} Suppose that $\alpha<\beta\in \omega_1$ and let $l=\rho(\alpha,\beta)$. If $F\in \mathcal{F}_{l}$ is such that $\{\alpha,\beta\}\subseteq F$, then $$\llbracket \alpha,\beta\rrbracket= \min(\,F_{\Xi_\beta(l)}\backslash R(F)\,).$$
\end{rem}

\begin{rem}\label{remarkbracketcaptured}As a consequence of the previous remark we have that if $D_0,D_1\in \text{Fin}(\omega_1)$ are disjoint sets for which $\{D_0,D_1\}$ is captured, then $$\llbracket D_0(i),D_1(i)\rrbracket=\llbracket D_0(j),D_1(j)\rrbracket$$
for all $i,j<|D_0|.$
\end{rem}
\begin{definition}[The square-bracket filter]Let $\mathcal{F}$ be a construction scheme. We define $\mathcal{V}(\mathcal{F})$ as the set of all $A\subseteq \omega_1$ for which there is $S\in [\omega_1]^{\omega_1}$  such that $$\llparenthesis S\rrparenthesis \subseteq A.$$
\end{definition}
Recall that a filter over $\omega_1$ is uniform whenever all of its elements are uncountable.
\begin{proposition}Let $\mathcal{F}$ be a $2$-capturing construction scheme. Then $\mathcal{V}(\mathcal{F})$ is a uniform filter.
\begin{proof}Since $\mathcal{F}$ is $2$-capturing, each element of $\mathcal{V}(\mathcal{F})$ is non-empty. In order to show that $\mathcal{V}(\mathcal{F})$ is a filter, it is enough to show that the family $\{\llparenthesis S\rrparenthesis\,:\,S\in[\omega_1]^{\omega_1}\}$ is downwards directed. Let $S,S'\in [\omega_1]^{\omega_1}$. We will prove that there is $A\in [\omega_1]^{\omega_1}$ so that $\llparenthesis A\rrparenthesis\subseteq \llparenthesis S\rrparenthesis \cap \llparenthesis S'\rrparenthesis.$ Just as in the Proposition \ref{lemmafilteruf}, we can consider a sequence $\langle \alpha_\xi,\beta_\xi\rangle_{\xi\in\omega_1}\subseteq S\times S'$ so that $\alpha_\xi<\beta_\xi<\alpha_{\zeta}$ for each $\xi,\zeta\in\omega_1$. We may assume without loss of generality that there is $k\in\omega$ so that the following properties hold for any two distinct $\xi<\mu\in\omega_1:$
\begin{enumerate}[label=$(\arabic*)$]
    \item $\rho(\alpha_\xi,\beta_\xi)<k$,
    \item $\lVert \alpha_\xi\rVert_k=\lVert \alpha_\mu \rVert_k$ and $\lVert \beta_\xi\rVert_k=\lVert \beta_\mu\rVert_k$. In particular, this means that \break$\rho(\beta_\xi,\beta_\mu)>k.$
\end{enumerate}
 We claim that $A=\{\beta_\xi\,:\,\xi\in\omega_1\}$ works. Trivially $\llparenthesis A\rrparenthesis \subseteq \llparenthesis S'\rrparenthesis$  because $A\subseteq S'$. In order to prove that $\llparenthesis A\rrparenthesis \subseteq \llparenthesis S\rrparenthesis$, let $\delta\in \llparenthesis A\rrparenthesis.$ Then there are $\xi<\mu\in\omega_1$ for which $\{\beta_\xi,\beta_\mu\}$ is captured at some level $l$ and $\llbracket \beta_\xi,\beta_\mu\rrbracket=\delta$. The properties (1) and (2) stated above imply that the family $\{\{\alpha_\xi,\beta_\xi\},\{\alpha_\mu,\beta_\mu\}\}$ is also captured at level $l$ (see Lemma \ref{capturedfamiliestosetslemma}). In particular, $\{\alpha_\xi,\alpha_\mu\}$ is captured. Furthermore, $\delta=\llbracket \beta_\xi,\beta_\mu\rrbracket=\llbracket \alpha_\xi,\alpha_\mu\rrbracket$ by means of the Remark \ref{remarkbracketcaptured}. As $\alpha_\xi,\alpha_\mu \in S$, we have shown that $\delta\in \llparenthesis S\rrparenthesis.$ 
 The only thing remaining to prove is that $\mathcal{V}(\mathcal{F})$ is uniform. This is a direct consequence of the fact that $\alpha\leq \llbracket \alpha,\beta\rrbracket \leq \beta$ for any two $\alpha<\beta\in \omega_1.$
\end{proof}    
\end{proposition}

\begin{theorem}Let $\mathcal{F}$ be a $2$ capturing construction scheme for which $\mathfrak{m}_\mathcal{F}>\omega_1$. Then $\mathcal{V}(\mathcal{F})$ is an ultrafilter.  
\begin{proof}Let $A\in \mathcal{V}(\mathcal{F})^+$. We will show that $A\in \mathcal{V}(\mathcal{F})$. For this purpose, let us consider the forcing $$\mathbb{P}=\{\,p\in[\omega_1]^{<\omega}\,:\,\llparenthesis p\rrparenthesis\subseteq A\,\}$$
ordered by reverse inclusion. It is straightforward that if $G\subseteq \mathbb{P}$ is an uncountable filter, then $S_G=\bigcup G$ is an uncountable subset of $\omega_1$ for which $\llparenthesis S_G\rrparenthesis \subseteq A$. In order to guarantee the existence of such filter, it is enough to prove that $\mathbb{P}$ is $ccc$ and $2$-preserves $\mathcal{F}$. This is because we are assuming that $\mathfrak{m}_\mathcal{F}>\omega_1$.\\

\noindent
\underline{Claim}: $\mathbb{P}$ is $ccc$ and $2$-preserves $\mathcal{F}$.
\begin{claimproof}[Proof of claim]Let $\mathcal{A}\in [\mathbb{P}]^{\omega_1}$ and $\nu:\mathcal{A}\longrightarrow \omega_1$ be an injective function. For each $p\in \mathcal{A}$ put $D_p=p\cup\{\nu(p)\}$$\alpha_p=\max(D_p)$. By refining $\mathcal{A}$ if necessary, we may assume that there are $j\in\omega$, $a<m_j$, $C\subseteq a+1$ and $b\in C$ so that the following conditions hold for any $p\in \mathcal{A}$:
\begin{enumerate}[label=$(\alph*)$]
    \item $\rho^{D_p}\leq j$,
    \item $\lVert \alpha_p\rVert_j=a, $
    \item $(\alpha_p)_j[C]=p$ and $(\alpha_p)_j(b)=\nu(p)$.
\end{enumerate}
Since $A\in \mathcal{V}(\mathcal{F})^+$, there are $p_0,p_1\in \mathcal{A}$ for which $\{\alpha_{p_0},\alpha_{p_1}\}$ is captured at some level $l$ and $\llbracket \alpha_{p_0},\alpha_{p_1}\rrbracket\in A$. According to Lemma \ref{capturedfamiliestosetslemma}, $\{D_{p_0},D_{p_1}\}$ is also captured at level $l.$ From this fact and by point (c), the same holds for $\{p_0,p_1\}$ and $\{\nu(p_0),\nu(p_1)\}$. In order to finish the claim, we will show that $q=p_0\cup p_1\in \mathbb{P}$. Indeed, let $\alpha<\gamma\in q$ be such that $\{\alpha,\gamma\}$ is captured.  If $\rho(\alpha,\gamma)<l$, then there is $i\in 2$ for which $\alpha,\gamma\in p_i$. Thus, $\llbracket \alpha,\gamma\rrbracket \in A$. On the other hand, if $\rho(\alpha,\gamma)=l$, let $F\in \mathcal{F}_l$ be such that $q\subseteq F$. By Remark \ref{remarksqbracket1}, we have that $\llbracket \alpha,\gamma\rrbracket= F_1\backslash R(F)=\llbracket \alpha_{p_0},\alpha_{p_1}\rrbracket\in A.$
\end{claimproof}

\end{proof}    
\end{theorem}
Suppose that $\mathcal{F}$ is as in the theorem above. Since $\omega_1$ is not measurable, there is a surjective function $\pi:\omega_1\longrightarrow \omega$ so that $\pi^{-1}[\{n\}]\not\in \mathcal{V}(\mathcal{F})$
for each $n\in\omega$. From this, it follows that the family $$\pi \mathcal{V}(\mathcal{F})=\{A\in \omega\,:\,\pi^{-1}[A]\in \mathcal{V}(\mathcal{F})\}$$
is a non-principal ultrafilter over $\omega$. We will now prove that it is also a Ramsey ultrafilter.
\begin{theorem}\label{pivframseytheorem}Let $\mathcal{F}$ be a $2$-capturing construction scheme for which $\mathfrak{m}_\mathcal{F}>\omega_1$. then $\pi \mathcal{V}(\mathcal{F})$ is a Ramsey ultrafilter.
\begin{proof}Let $\langle P_i\rangle_{i\in\omega}$ be a partition of $\omega$ so that $P_i\not\in \pi\mathcal{V}(\mathcal{F}) $ for each $\in \omega$. We define the forcing $\mathbb{P}$ as the set of all $p\in [\omega_1]^{\omega}$ so that $$|\pi[ \llparenthesis p\rrparenthesis]\cap P_i|\leq 1$$
for each $i\in\omega$. We order $\mathbb{P}$ with the reverse inclusion.
Note that if $G\subseteq\mathbb{P}$ is an uncountable filter, then $S_G$ is an uncountable subset of $\omega_1$ such that $|\pi[\llparenthesis S_G\rrparenthesis]\cap P_i|\leq 1$ for each $i$. The existence of such filter follows from the following claim.\\\\
\noindent\underline{Claim}: $\mathbb{P}$ is $ccc$ and $2$-preserves $\mathcal{F}$.
\begin{claimproof}[Proof of claim]
Let $\mathcal{A}$ be an uncountable subset of $\mathbb{P}$ and $\nu:\mathcal{A}\longrightarrow \omega_1$ be an injective function. Given $p\in \mathcal{A}$, let $D_p=p\cup\{\nu(p)\}$ and $\alpha_p=\max(p)$.  By refining $\mathcal{A}$ if necessary, we may assume that there are $j\in\omega$, $a<m_j$, $C\subseteq a+1$ and $b\in C$ so that the following conditions hold for any $p\in \mathcal{A}$:
\begin{enumerate}[label=$(\alph*)$]
    \item $\rho^{D_p}\leq j$,
    \item $\lVert \alpha_p\rVert_j=a, $
    \item $(\alpha_p)_j[C]=p$ and $(\alpha_p)_j(b)=\nu(p)$.
\end{enumerate}
Furthermore, we may suppose that there is $E\subseteq \omega$ so that $\pi[\llparenthesis p\rrparenthesis]=E$ for each $p\in \mathcal{A}$. Let $I=\{i\in\omega\,:\,E\cap P_i\not=\emptyset\}$. Note that $I$ is finite. Therefore, $$A=\omega\backslash (\bigcup\limits_{i\in I}P_i)\in \pi \mathcal{V}(\mathcal{F}).$$
In particular, $A\cap \pi[ \llparenthesis \alpha_p\,:\,p\in \mathcal{A}\rrparenthesis]\not=\emptyset$. Thus, there are $p_0,p_1\in \mathcal{A}$ so that $\{\alpha_{p_0},\alpha_{p_1}\}$ is captured at some level $l$ and $\pi ( \llbracket \alpha_{p_0}\alpha_{p_1}\rrbracket )\in A$. By arguing in the same way as in the previous theorem, we can conclude that both $\{\nu(p_0),\nu(p_1)\}$ and $\{p_0,p_1\}$ are captured. Even more, for $q=p_0\cup p_1$ we have that $$\llparenthesis q\rrparenthesis=\{\llbracket \alpha_{p_0},\alpha_{p_1}\rrbracket\}\cup \llparenthesis p_0\rrparenthesis \cup \llparenthesis p_1\rrparenthesis.$$
In virtue of this,  $\pi[\llparenthesis q\rrparenthesis]=\{\pi(\llbracket \alpha_{p_0},\alpha_{p_1}\rrbracket)\}\cup E$. As $\pi[\llbracket \alpha_{p_0},\alpha_{p_1}\rrbracket]\not\in A$, it easily follows that $q\in \mathbb{P}$. This finishes the proof.
\end{claimproof}
\end{proof}
\end{theorem}

\section{P-ideal Dichotomy and capturing axioms}\label{sectionpidschemes}

Let $\mathcal{I}$ be an ideal over a set $X$. We say that $\mathcal{I}$ is a \textit{$P$-ideal} if for each countable $\mathcal{A}\subseteq \mathcal{I}$ there is $B\in \mathcal{I}$ such that $A\subseteq^*B$ for every $A\in \mathcal{A}$.
Any such set $B$ is called a \textit{pseudo-union} of $\mathcal{A}$. The $P$-ideal dichotomy (PID) is a well-known consequence of PFA regarding $P$-ideals. The third author introduced it in its current form in \cite{dichotomypidfirst}. There he showed that PID is consistent with CH, although PFA is not. A less general version of PID was previously formulated by Abraham and  the third author in \cite{partitionpropertiesch}.

   \begin{definition}Let $\mathcal{I}$ be an ideal over a set $X$. We define the \textit{orthogonal ideal of $\mathcal{I}$}, namely  $\mathcal{I}^\perp$,  as the set of all $A\subseteq X$ such that $A\cap B=^*\emptyset$ for any $B\in \mathcal{I}.$  
\end{definition}

\noindent
{\textbf{$P$-ideal Dichotomy [PID]:}} Let $X$ be a set and $\mathcal{I}$ be a $P$-ideal over $X$ such that $[X]^{<\omega}\subseteq \mathcal{I}\subseteq [X]^{\leq \omega}$. Then one of the two following conditions hold:\begin{itemize}
    \item There is $Y\in [X]^{\omega_1}$ so that $[Y]^\omega\subseteq \mathcal{I}$.
    \item There is $\{ Z_n\,:\,n\in\omega\}\subseteq\mathcal{I}^\perp$ for which $X=\bigcup\limits_{n\in\omega} Z_n.$
\end{itemize}

Many of the objects that can be constructed using capturing construction schemes already contradict to some degree the $P$-ideal dichotomy. For example, since the Suslin Hypothesis follows from PID (see \cite{partitionpropertiesch} and \cite{dichotomypidfirst}), this axiom is incompatible with FCA. In addition, since PID implies that any gap is indestructible, $CA_3$ also contradicts PID. However, non of the consequences of  $CA_2$ which are known so far are enough to contradict the $P$-ideal dichotomy. For now, the best we have is that $PID+\min(\mathfrak{b},\,cof(\mathcal{F}_\sigma)\,)>\omega_1$ is incompatible with the existence of a $2$-capturing construction scheme. This is because such assertion is equivalent to the nonexistence of a sixth Tukey type (see \cite{CombinatorialDichotomiesandCardinalInvariants}). The purpose of this section is to prove that PID is in fact incompatible with $CA_2$.

For the remainder of this section, we will assume that PID holds and that $\mathcal{F}$ is a construction scheme. Our plan is to define a $P$-ideal $\mathcal{I}$ from $\mathcal{F}$ that contradicts PID whenever $\mathcal{F}$ is $2$-capturing. This ideal $\mathcal{I}$ will be defined as family of countable sets of the orthogonal ideal of some other ideal $\mathcal{H}$. When applying the second alternative of PID to the ideal $\mathcal{I}$, we will need to calculate who is $\mathcal{I}^\perp$. In general, this may be a difficult task. However, the following concept will facilitate these calculations.

\begin{definition}Let $\mathcal{H}$ be an ideal of countable sets in $\omega_1$. We will say that $\mathcal{H}$ is \textit{Fr\'echet in $[\omega_1]^{\omega}$} if for any $A\in \mathcal{H}^+\cap[\omega_1]^\omega$, there is an infinite $B\subseteq A$ such that $B\in \mathcal{H}^\perp$.
    
\end{definition}

\begin{proposition}Let $\mathcal{H}$ be an ideal of countable sets in $\omega_1$ and assume that $\mathcal{H}$ is Fr\'echet in $[\omega_1]^{\omega_1}$. Then $(\mathcal{H}^\perp\cap [\omega_1]^{\leq\omega})^\perp\cap [\omega_1]^{\omega}=\mathcal{H}\cap [\omega_1]^{\omega}$.
\begin{proof}We will only show that the  inclusion from left to right holds. For this purpose, let $A\in (\mathcal{H}^\perp\cap[\omega_1]^{\leq \omega})^\perp\cap[\omega_1]^{\omega}$. Assume towards a contradiciton that $A\not\in \mathcal{H}$. Then $A\in \mathcal{H}^+\cap [\omega_1]^\omega$. Since $\mathcal{H}$ is Fr\'echet, there is an infinite $B\subseteq A$ such that $B\in \mathcal{H}^\perp$. In particular, $A\cap B$ is infinite, and this contradicts the fact that $A\in (\mathcal{H}^\perp)^\perp$. Thus, the proof is over.
\end{proof}
\end{proposition}
We now define the ideal $\mathcal{I}$ which will be used in the proof of Theorem \ref{fnotcapturingpidtheorem}.

\begin{definition} Given $n\in\omega$ and $\alpha\in \omega_1$, we define $$H_n(\alpha)=\{\xi<\alpha\,:\,\forall m>n\,(\,\Xi_\alpha(m)=-1\textit{ or }\Xi_\xi(m)\leq \Xi_\alpha(m)\,)\}.$$
Additionally, we define the ideal generated by the family $\{H_n(\alpha)\,:\,n\in\omega\textit{ and }\alpha\in\omega_1\}$ as $\mathcal{H}$. Lastly, we define $\mathcal{I}=\mathcal{H}^\perp\cap [\omega_1]^{\leq \omega}.$
\end{definition}
\begin{rem}$[\omega_1]^{<\omega}$ is included in both $\mathcal{H}$ and $\mathcal{I}.$
    
\end{rem}
The following lemma is easy.
\begin{lemma}\label{lemmaHmpideal1}Let $\alpha,\beta,\gamma<\omega_1$ and $m,n\in \omega$. Then:
\begin{enumerate}[label=$(\arabic*)$]
    \item If $m<n$, then $H_m(\alpha)\subseteq H_n(\alpha).$
    \item If $\beta<\gamma$ and $n>\rho(\beta,\gamma)$, then $H_n(\gamma)\cap \beta\subseteq H_n(\beta).$
\end{enumerate}
\begin{proof}The point (1) follows directly from the definition. In order to prove (2) we will show that $\omega_1\backslash H_n(\beta)\subseteq \omega_1 \backslash (H_n(\gamma)\cap \beta)$. Let $\xi \in \omega_1\backslash H_n(\beta) $. If $\xi\geq \beta$ we are done, so let us assume that $\xi<\beta.$ Then there is $m>n$ for which $\Xi_\xi(m)>\Xi_\beta(m)\geq 0$. Note that $m>\rho(\beta,\gamma)$. Thus, by the point (c) of Lemma \ref{lemmaxi} we conclude that $\Xi_\beta(m)=\Xi_\gamma(m)$. In this way, $\xi\not\in H_n(\gamma)$. This finishes the proof.
\end{proof}
    
\end{lemma}
\begin{rem}By the point (1) of the previous lemma we have that for any $A\in \mathcal{H}$ there are $n\in\omega$ and $\alpha_1,\dots,\alpha_n\in \omega_1$ for which $A\subseteq \bigcup\limits_{i<n}H_n(\alpha_i).$
    
\end{rem}

\begin{definition}Given $\beta\in \omega_1$ we define $\mathcal{H}|_\beta$ as the ideal generated by the family $\{H_n(\alpha)\,:\,\alpha\leq \beta \textit{ and }n\in\omega\}.$
    
\end{definition}
\begin{lemma}\label{lemmaHmpideal2}Let $\beta<\omega_1$ and $A\subseteq \beta$.
\begin{enumerate}[label=$(\alph*)$]
    \item If $A\in (\mathcal{H}|_\beta)^{\perp}$, then $A\in \mathcal{H}^\perp.$ Furthermore, since $A$ is countable then $A\in \mathcal{I}.$
    \item If $A\in (\mathcal{H}|_\beta)^{+}$, then $A\in \mathcal{H}^+.$
\end{enumerate}
\begin{proof}\begin{claimproof}[Proof of $(a)$] It is enough to show that $A\cap H_n(\gamma)=^*\emptyset$ for any $\beta\leq\gamma\in \omega_1$ and each $n\in\omega$. Furthermore, we may assume that $n>\rho(\beta,\gamma)$ due to the point (1) of Lemma \ref{lemmaHmpideal1}. Indeed, by virtue of the point (2) of the same lemma, we have that $$A\cap H_n(\gamma)=A\cap \beta\cap H_n(\gamma)\subseteq A\cap H_n(\beta)=^*\emptyset.$$
\end{claimproof}
\begin{claimproof}[Proof of $(b)$]Suppose towards a contradiction that $A\notin \mathcal{H}^+$. Then there are $n\in\omega$ and $\alpha_1,\dots,\alpha_n\in \omega_1$ for which $A\subseteq \bigcup\limits_{i<n}H_n(\alpha_i).$ Without any loss of generality we may assume that there is $j<n$ so that, for any $i<n$, $\alpha_i\leq \beta $ if and only if $i<j$. Furthermore, we may assume that $n>\rho(\beta,\alpha_i)$ for any $i\geq j$. From this, it follows that $$A\subseteq\big( \bigcup\limits_{i<j}H_n(\alpha_i)\big)\cup \big(\bigcup\limits_{j\leq i<n}H_n(\alpha_i)\cap \beta\big)\subseteq \big(\bigcup\limits_{i<j}H_n(\alpha_i)\big)\cup \big(\bigcup\limits_{j\leq i<n}H_n(\beta)\big). $$
In this way, $A\in \mathcal{H}|_\beta$ which is a contradiction.  
\end{claimproof} 
\end{proof}
    
\end{lemma}

\begin{proposition}\label{propIFrechet}$\mathcal{H}$ is Fr\'echet in $[\omega_1]^\omega.$
\begin{proof}Let $A\in \mathcal{H}^+\cap [\omega_1]^\omega$ and $\beta\in\omega_1$ be such that $A\subseteq \beta$. Note that $A\in (\mathcal{H}|_\beta)^+$. Since $\mathcal{H}|_\beta$ is a countably generated ideal, there is $B\in [A]^\omega$ such that $B\in (\mathcal{H}|_\beta)^\perp$. According to the point $(a)$ of Lemma \ref{lemmaHmpideal2}, $A\in \mathcal{I}$. Thus, the proof is over.
    
\end{proof}
\end{proposition}

\begin{proposition}\label{propIpideal}
$\mathcal{I}$ is a $P$-ideal.
\begin{proof}Let $\langle A_n\rangle_{n\in\omega}$ be an increasing sequence of elements of $\mathcal{I}$. Then there is $\beta\in\omega_1$ for which $\bigcup\limits_{n\in\omega}A_n\subseteq \beta$. Let us enumerate $\beta+1$ as $\langle \beta_n\rangle_{n\in\omega}$. We then define $$A=\bigcup\limits_{n\in\omega}\big(A_n\backslash \big(\bigcup\limits_{i\leq n}H_n(\beta_i)\big)\big).$$
By the Lemma \ref{lemmaHmpideal1}, it follows that $A$ is a pseudo-union of $\langle A_n\rangle_{n\in\omega}$ and $A\in (\mathcal{H}|_\beta)^\perp$. In particular, this means that $A\in \mathcal{I}$ due to the point (a) of Lemma \ref{lemmaHmpideal2}.
\end{proof}
\end{proposition}
\begin{theorem}[Under PID]\label{fnotcapturingpidtheorem}$\mathcal{F}$ is not $2$-capturing.
\begin{proof} Suppose towards a contradiction that $\mathcal{F}$ is $2$-capturing. According to the Proposition \ref{propIpideal}, $\mathcal{I}$ is a $P$-ideal. Since we are assuming PID, there are two possibilities for $\mathcal{I}$. We will finish the proof by showing that both of them lead to a contradiciton.\\\\
\noindent
\underline{Case 1}: There is an uncountable $Y\subseteq \omega_1$ such that $[Y]^{\omega}\subseteq \mathcal{I}.$
\begin{claimproof}[Non satisfaction of case]Consider the coloring $c:[Y]^2\longrightarrow 2$ given by:
$$c(\{\alpha,\beta\})=0\text{ if and only if }\{\alpha,\beta\}\text{ is not captured}.$$
Using the partition relation $\omega_1\longrightarrow (\omega_1,\omega+1)^2_2$ and the fact that $\mathcal{F}$ is $2$-capturing, we can conclude that there is an increasing sequence $\langle \alpha_i\rangle_{i<\omega+1}\subseteq Y$ such that $\{\alpha_i,\alpha_j\}$ is captured for any two distinct $i,j< \omega+1$. In particular, $\{\alpha_i,\alpha_\omega\}$ is captured for any $i<\omega$. By Lemma \ref{lemmaxi}, we get that $A=\{\alpha_i\,:\,i<\omega\}\subseteq H_0(\alpha_\omega)$. In other words, $A\in \mathcal{I}$. This is a contradiction, so this case can not occur.
\end{claimproof}
\noindent
\underline{Case 2}: There is $Y$ an uncountable subset of $\omega_1$ such that $[Y]^{\omega}\subseteq \mathcal{I}^\perp.$
\begin{claimproof}[Non satisfaction of case]
According to Proposition \ref{propIFrechet}, $\mathcal{H}$ is Fr\'echet over $[\omega_1]^\omega$. From this it follows that  $[Y]^{\omega}\subseteq \mathcal{H}$. Let $$P=\{l\in\omega\,:\,\exists \alpha\in Y\,(\,\Xi_\alpha(l)>0\,)\}.$$ By refining $Y$,  we may assume that for any $l\in P$ there are uncountable many $\alpha\in Y$ for which $\Xi_\alpha(l)>0$. In this way, we can recursively construct sequences $\langle A_\beta\rangle_{\beta\in\omega_1}\subseteq [Y]^\omega$, $\langle B_\beta\rangle_{\beta\in \omega_1}\subseteq \text{Fin}(\omega_1)$ and $\langle n_\beta\rangle_{\beta\in\omega_1}$ so that the following properties hold for any two $\beta<\gamma\in \omega_1$:
\begin{enumerate}[label=$(\arabic*)$]
    \item $A_\beta<B_\beta<A_\gamma$. In particular $B_\beta\cap B_\gamma=\emptyset.$
    \item For all $l\in P$ there is $\alpha\in A_\beta$ for which $\Xi_\alpha(l)>0$.
    \item $A_\beta\subseteq \bigcup\limits_{\alpha\in B_\beta}H_{n_\beta}(\alpha).$
    \item $B_\beta\cap Y\not=\emptyset.$ 
\end{enumerate}
By refining such sequence, we can suppose that there is $n\in\omega$ such that $n_\beta=n$ for any $\beta\in\omega_1$. Since $\mathcal{F}$ is $2$-capturing, there are $\beta<\gamma\in\omega_1$ such that $\{B_\beta,B_\gamma\}$ is captured at some level $l>n$. Now, $B_\beta$ and $B_\gamma$ are disjoint. Thus, $\Xi_\alpha(l)=1$ for any $\alpha\in B_\gamma$. Therefore, $l\in P$ due to the point (4). By virtue of the point (2), there must be $\xi\in A_\beta$ such that $\Xi_\xi(l)>0$. On the other hand, $\Xi_\alpha(l)=0$ for all $\alpha\in B_\beta.$ Consequently, $\xi \not\in H_n(\alpha)$ for any such $\alpha$. This contradicts the point (3), so this case can not occur.  
\end{claimproof}  
\end{proof}
\end{theorem}
\begin{corollary}(Under PID)\label{coropidnotca2} There are no $2$-capturing construction schemes.
    
\end{corollary}

\section{Open problems}\label{sectionproblems}
The main objects studied through out this text were the WLI gaps. We wonder how much PFA is needed to imply that every gap is of this form. A particular instance of this line of  research would be the following:
\begin{problem}Does PID imply that every gap is weakly levelwise-separable? If not,  what about PID+MA?
    
\end{problem}
In this paper, we solely dealt  with gaps of size $\omega_1$. However, the notion of  levelwise-inseparabilty makes sense for gaps of other sizes. Having this generalization in mind we ask:
\begin{problem}For which pairs of ordinals $(\kappa,\lambda)\not=(\omega_1,\omega_1)$ is consistent that there is a $(\kappa,\lambda)$-gap which is not weakly levelwise-inseparable? For which of those pairs does the statement \say{There is a $(\kappa,\lambda)$-gap and every $(\kappa,\lambda)$-gap is weakly levelwise-inseparable} is consistent?
\end{problem}
 We know that Suslin trees exist if a fully capturing schemes (of a certain type) exists, and that $2$-capturing schemes are not powerful enough to imply the existence of these objects. However, we still do not know what is the situation for $n$-capturing (or capturing) schemes. Compare this with the fact that destructible gaps do exist if we assume that there is a $3$-capturing scheme.
\begin{problem}Let $n>2$. Does the inequality $\mathfrak{m}^n_\mathcal{F}>\omega_1$ imply that there are no Suslin trees?
\end{problem}

A theorem of Baumgartner \cite{basestrees} states that it is consistent (under large cardinal hypotheses) that every Aronszajn tree has a base of cardinality $\omega_1$. This in turn implies that consistently the filter $\mathcal{U}(T)$ associated to a coherent tree has a (filter) base of size $\omega_1.$  We wonder whether the same is true for the filter $\mathcal{V}(\mathcal{F})$ defined in Section \ref{sectionramseyultafilters}.
\begin{problem}Is it consistent that $\mathcal{V}(\mathcal{F})$ is a filter is generated by $\omega_1$-many sets?
A positive answer to the question above may shed some light to Kunen's open problem of whether there can be an ultrafilter $\mathcal{U}$ over $\omega_1$ whose character is less than $2^{\omega_1}$.
    
\end{problem}

   \bibliographystyle{plain}
\bibliography{phDThesis}
\end{document}